\documentclass[11pt,a4paper]{article}
\usepackage{amsfonts}
\usepackage{amsmath}
\usepackage{amssymb}
\usepackage{graphicx}
\usepackage{tikz}
\usetikzlibrary{snakes}

\newtheorem{theorem}{\indent Theorem}

\newtheorem{lemma}{\indent Lemma}

\newtheorem{proposition}{\indent Proposition}
\newtheorem{remark}{Remark}

\newcommand{\ud}{\mathrm{d}}

\def\kraj{\hfill\rule{6pt}{6pt}}
\def\diag{\mathop{\rm diag}}
\def\deg{\mathop{\rm deg}}
\def\rank{\mathop{\rm rank}}
\def\rk{\mathop{\rm rk}}
\def\sgn{\mathop{\rm sgn}}
\def\lcm{\mathop{\rm lcm}}
\def\tr{\mathop{\rm tr}}
\def\gcd{\mathop{\rm gcd}}
\def\Hom{\mathop{\rm Hom}}
\def\Id{\mathop{\rm Id}}
\def\res{\mathop{\rm res}}
\def\F{\mathbb{F}}
\def\R{\mathbb{R}}
\def\Q{\mathbb{Q}}
\def\Z{\mathbb{Z}}
\def\C{\mathbb{C}}
\def\S{\mathbb{S}}
\def\ka{\mathbb{K}}
\def\E{\mathcal{E}}
\def\U{\mathcal{U}}
\def\B{\mathcal{B}}
\def\Ker{\mathop{\rm Ker}}
\def\max{\mathop{\rm max}}
\def\Seq{\mathbf{\rm Seq}}
\def\Bot{\mathop{\rm Bot}}
\def\Top{\mathop{\rm Top}}
\def\dun{\dot{\U}^+_n}
\def\un{{\U}^+_n}
\def\K{\mathbb{K}}
\def\N{\mathbb{N}}
\def\u3{U^+_q(\mathfrak{sl}_3)}
\def\sl{\mathfrak{sl}}
\author{Marko Sto\v si\'c}
\title{Indecomposable 1-morphisms of $\dot{\mathcal{U}}^+_3$ and the canonical basis of ${U}_q^+(\mathfrak{sl}_3)$}

\begin{document}
\maketitle

\begin{abstract}
We compute the indecomposable objects of $\dot{\mathcal{U}}^+_3$ -- the category that categorifies the positive half of the quantum $\mathfrak{sl}_3$, and we decompose an arbitrary object into indecomposable ones. 
On decategorified level we obtain the Lusztig's canonical basis of the positive half ${U}^+_q(\mathfrak{sl}_3)$ of the quantum $\mathfrak{sl}_3$. We also categorify the higher quantum Serre relations 
in ${U}_q^+(\mathfrak{sl}_3)$, by defining a certain complex in the homotopy category of $\dot{\mathcal{U}}^+_3$ that is homotopic to zero. We work with the category $\dot{\mathcal{U}}^+_3$ that is defined over the ring of integers. This paper is based on the (extended) diagrammatic calculus introduced to categorify quantum groups. 
\end{abstract}

\section{Introduction}\label{sec1}

In recent years there has been a lot of work on a diagrammatic categorification of quantum groups, initiated by Lauda's diagrammatic categorification \cite{sl2} of the Lusztig's idempotented version of $\dot{U}_q(\sl_2)$. This was extended by Khovanov and Lauda in \cite{kl3} to $\dot{U}_q(\sl_n)$ and also in \cite{kl2} to the positive half of an arbitrary quantum group $U^{+}_q(\mathfrak{g})$. In \cite{w1,w2}, Webster modified the construction, and obtained a  categorification of an arbitrary quantum group.\\

General framework of these constructions is to define a certain $2$-category $\U$ whose $1$-morphisms categorify the generators of a quantum group, and whose $2$-morphisms are $\ka$-linear combinations of certain planar diagrams modulo local relations, with $\ka$ being a field. Then a $2$-category $\dot{\U}$ is defined as the Karoubi envelope of the $2$-category $\U$, i.e. as the smallest $2$-category containing ${\U}$ in which all idempotents ($2$-morphisms) split. Finally, it is shown that the split Grothendieck group of $\dot{\U}$ is isomorphic to the corresponding quantum group.

In the case of the categorification of the positive half of quantum groups the  $2$-categories $\U$ and $\dot{\U}$ have a single 
object. Thus one can see them as monoidal $1$-categories. Since in this paper we are interested in categorifications of positive halves of quantum groups, we shall always assume that $\U$ and $\dot{\U}$ are monoidal ($1$-)categories.\\

In order to be able to use diagrammatical calculus in $\dot{\U}$ directly  (and not just in $\U$), in \cite{thick} the extension of the calculus -- so called thick calculus -- was introduced in the case of quantum $\mathfrak{sl}_2$. The lines labelled $a$ correspond to objects of $\dot{\U}$ that categorify the divided powers $E^{(a)}$ of  the generators of the quantum group. The consequence of \cite{thick}, and of the thick calculus, is that now one can take $\Z$-linear combinations of planar diagrams as morphisms of $\U$.\\

In this paper, we use thick calculus to study the properties of the category $\dot{\U}^{+}_n$ that categorifies the positive half of the quantum $\sl_n$. The thick calculus can be extended directly to include the categorification of $\sl_n$ (see \cite{thickn}). In particular, the category $\dot{\U}^{+}_n$ is defined over the ring of integers.\\

In Section \ref{sec4} of this paper, we compute the indecomposable objects of $\dot{\U}^{+}_3$. More precisely, in Theorem \ref{te3}, we show that these are the following:
$$\mathcal{B}=\left\{ \E_1^{(a)}\E_2^{(b)}\E_1^{(c)}\{t\}, \E_2^{(a)}\E_1^{(b)}\E_2^{(c)}\{t\} \mid b\ge a+c,\quad a, c\ge 0,\,\,\,t\in\Z\right\},$$
with $\E_1^{(a)}\E_2^{(a+c)}\E_1^{(c)}\{t\}\cong \E_2^{(c)}\E_1^{(a+c)}\E_2^{(a)}\{t\}$, for $a,c\!\ge\! 0$, $t\in\Z$. Moreover, we prove that an arbitrary object of $\dot{\U}^{+}_3$ can be decomposed as a direct sum of the elements of $\mathcal{B}$. 
%In particular, the composition of 
%two $1$-morphisms from $\mathcal{B}$ can be decomposed as a direct sum of the elements of $\mathcal{B}$.

The main result (Theorem \ref{gldect}), is the categorification of the  ${U}^{+}_q(\sl_3)$ relation:
\begin{equation}E_1^{(a)}E_2^{(b)}E_1^{(c)}=\sum_{{\scriptsize{\begin{array}{c}p+r=b\\
p\le c\\r\le a\end{array}}}}\left[\begin{array}{c}a+c-b\\c-p\end{array}\right]E_2^{(p)}E_1^{(a+c)}E_2^{(r)},\quad \textrm{for } b\le a+c.\label{intro}\end{equation}

By decategorifying the set of indecomposables from $\mathcal{B}$ with no shifts, we obtain the set 
$$B=\left\{ E_1^{(a)}E_2^{(b)}E_1^{(c)},\,\,\, E_2^{(a)}E_1^{(b)}E_2^{(c)} \mid b\ge a+c,\quad  a, c\ge 0\right\},$$
with $E_1^{(a)}E_2^{(a+c)}E_1^{(c)}= E_2^{(c)}E_1^{(a+c)}E_2^{(a)}$, for $a,c\!\ge\! 0$.
The set $B$ is the Lusztig's canonical basis of ${U}^{+}_q(\sl_3)$ (see \cite{luszr}), and one of its remarkable properties is that its structure constants are from $\N[q,q^{-1}]$.\\

Thus, in this way, we have proved that the indecomposable objects of $\dot{\U}^{+}_3$ lift the canonical basis of ${U}^{+}_q(\sl_3)$. We note once again that we are working in the category that is defined over the ring of integers $\Z$, i.e. the $1$-morphisms are $\Z$-linear combinations of certain planar diagrams. 

Previous results on this topic were obtained in the case when the category is defined over a field, i.e. when the $1$-morphisms are $\ka$-linear combinations of planar diagrams, for some field $\ka$. The result that the indecomposable objects of that category lift the Lusztig canonical basis in the case of $\sl_3$ was obtained by Khovanov and Lauda \cite{kl1}. Furthermore, Brundan and Kleshchev \cite{bk} have extended the result to the case of affine $\sl_n$. Finally, for $\ka=\C$, Varagnolo and Vasserot \cite{vv} proved this fact for any simply-laced $\mathfrak g$.\\

The other goal of this paper is to categorify the higher quantum Serre relations for $E_1$ and $E_2$. The higher quantum Serre relations for the generators $E_1$ and $E_2$ (and also for the generators $E_r$ and $E_s$ of the positive half of an arbitrary quantum group $U^{+}_q(\mathfrak{g})$ with $r\cdot s=-1$) are 
\begin{equation}
\sum_{i=0}^m (-1)^i q^{\pm(m-n-1)i}E_1^{(m-i)}E_2^{(n)}E_1^{(i)}=0,\textrm{ for } m>n>0.\label{intro2}
\end{equation}

The relation (\ref{intro2}) can be obtained by summing appropriately some relations of the form (\ref{intro}). Thus from the categorification of (\ref{intro}) (Theorem \ref{gldect}), one can obtain a decomposition that lifts (\ref{intro2}).

However, in Section \ref{sec5}, we give a direct and simple categorification of the higher quantum Serre relation in the homotopy category of $\dot{\U}^{+}_3$ -- the category of complexes in $\dot{\U}^{+}_3$, modulo homotopies. 

Since the higher quantum Serre relations have the form of an alternating sum, it is natural to look for a categorification in the form of a complex of objects of $\dot{\U}^{+}_3$ that lift the summands of (\ref{intro2}). In Theorem \ref{te5}, we define such a complex and show that it is homotopic to zero. Moreover, the differentials and the homotopies have a particularly simple form.\\

Finally, all results from this paper about the generators $E_1$ and $E_2$, and objects $\E_1$ and $\E_2$, are valid in an arbitrary quantum group $U_q(\mathfrak{g})$ and in the categorification of its positive half, for the generators $E_r$ and $E_s$, and for $\E_r$ and $\E_s$, respectively, with $r$ and $s$ satisfying $r\cdot s=-1$.

\section{${U}^+_q(\mathfrak{sl}_n)$}\label{sec2}

In this section we give basic definitions of the positive half of quantum 
$\mathfrak{sl}_n$ -- denoted ${U}^+_q(\mathfrak{sl}_n)$. We also give some of its combinatorial properties  in the case $n=3$. 
These properties are also valid for any two generators $E_i$ and $E_j$ of an arbitrary quantum group when $i\cdot j=-1$.\\

Let $n\ge 2$ be fixed. The index set of  quantum $\mathfrak{sl}_n$ is $I=\{1,2,\ldots,n-1\}$. 
The inner product is defined on $\Z[I]$ by setting that for $i,j\in I$:
$$ i\cdot j= \left\{\begin{array}{rl} 2,&\quad i=j\\ -1,&\quad |i-j|=1\\
0,&\quad |i-j|\ge 2\end{array}\right.$$

${U}^+_q(\mathfrak{sl}_n)$ is a  $\Q(q)$-algebra generated by $E_1,E_2,\ldots,E_{n-1}$ modulo relations:
\begin{eqnarray}
E_i^2E_j+E_jE_i^2&=&[2]E_iE_jE_i,\quad i\cdot j=-1,\\
E_iE_j&=&E_jE_i, \quad\quad i\cdot j =0.
\end{eqnarray}

The divided powers of the generators are defined by 
$$E_i^{(a)}:=\frac{E_i^a}{[a]!}, \quad a\ge 0,\,\,\,\, i=1,\ldots,n-1.$$

The divided powers satisfy:
\begin{eqnarray}
E_i^{(a)}E_j^{(b)}&=&E_j^{(b)}E_i^{(a)},\quad i\cdot j=0,\label{ko}\\
E_{i}^{(a)}E_{i}^{(b)}&=&\left[\begin{array}{c}a+b\\a\end{array}\right]E_{i}^{(a+b)},\label{eiab}
\end{eqnarray}
and  the quantum Serre relations
\begin{equation}
E_i^{(2)}E_j+E_jE_i^{(2)} =  E_iE_jE_i,\quad i\cdot j = -1.\label{kser}
\end{equation}
\vskip 0.3cm
The integral form ${}_{\Z}{U}^+_q(\mathfrak{sl}_n)$ is the 
$\Z[q,q^{-1}]$-subalgebra of ${U}^+_q(\mathfrak{sl}_n)$ generated by $E_{i}^{(a)}$, for all $i=1,\ldots,n-1$ and $a\ge 0$.

\subsection{Combinatorics of $U^+_q(\mathfrak{sl}_3)$}

The higher quantum Serre relations (Chapter 7 of \cite{luszb}) for $E_1$ and $E_2$ are:
\begin{equation}\label{hqs}
\sum_{r=0}^m {(-1)^r q^{\pm (m-n-1)r}\,\,\, E_1^{(m-r)}E_2^{(n)}E_1^{(r)}}=0, \quad m>n>0,
\end{equation}
and analogously with $E_1$ and $E_2$ switched.
In particular, the quantum Serre relations are obtained for $m=2$ and $n=1$.\\

\begin{proposition}[\cite{luszb}, Lemma 42.1.2.(d)]\label{prgl}
For any three nonnegative integers $a,b,c$, with $b\le a+c$, we have:
\begin{equation}\label{glavna}
E_1^{(a)}E_2^{(b)}E_1^{(c)}=\sum_{{\scriptsize{\begin{array}{c} p+r=b \\ p\le c \\ r \le a\end{array}}}}
\left[\begin{array}{c}a+c-b\\ c-p \end{array}\right]  E_2^{(p)}E_1^{(a+c)}E_2^{(r)}. 
\end{equation}
In particular, for $b=a+c$, we have:
\begin{equation}
E_1^{(a)}E_2^{(a+c)}E_1^{(c)}=E_2^{(c)}E_1^{(a+c)}E_2^{(a)}.
\end{equation}
Finally, both formulas are valid when $E_1$ and $E_2$ interchange places.
\end{proposition}

We note that the higher quantum Serre relations (\ref{hqs}) follow from (\ref{glavna}), together  
with the following well-known relation of quantum binomial coefficients, valid for any nonnegative integer $N$:
\begin{equation}
\sum_{k=0}^N {(-1)^k q^{\pm (N-1)k} \left[\!\!\!\begin{array}{c}N \\k \end{array}\!\!\!\right]}=0.
\end{equation}

All relations from above for $E_1$ and $E_2$ are also valid in an 
arbitrary quantum group for two generators $E_i$ and $E_j$, with $i\cdot j=-1$.

\subsection{Monomials in $U^+_q(\mathfrak{sl}_3)$}
By a monomial, we mean a vector of the form $E_1^{(a_1)}E_2^{(b_2)}E_1^{(a_2)}\ldots E_1^{(a_n)}E_2^{(b_n)}$. The number of nonzero exponents $a_i$ and $b_i$ we call the length of the vector. Let $B$ be the following set of monomials:
\begin{equation}
B=\{E_1^{(a)}E_2^{(b)}E_1^{(c)},\,\,\,E_2^{(a)}E_1^{(b)}E_2^{(c)}| \quad b\ge a+c,\,\,\, a,b,c\ge 0\}, 
\end{equation}
where we have $E_1^{(a)}E_2^{(a+c)}E_1^{(c)}=E_2^{(c)}E_1^{(a+c)}E_2^{(a)}$, for all $a,c\ge 0$. 

The set $B$ is the Lusztig's canonical basis of $\u3$, and its structure constants are from $\N[q,q^{-1}]$. 

We also have the following:

\begin{theorem}\label{mon}
Every monomial from $U^+_q(\mathfrak{sl}_3)$ can be expressed as a linear combination of  vectors from $B$, with coefficients from $\N[q,q^{-1}]$.
\end{theorem} 

\textbf{Proof:}

By induction on length. By Proposition \ref{prgl} any monomial of length at most 3 can be expressed as a linear combination of the vectors from $B$, which proves the base of induction.

Now, suppose that a monomial $v$ has length $l$, with $l\ge 4$. Then it contains a piece of the form 
$E_1^{(a)}E_2^{(b)}E_1^{(c)}E_2^{(d)}$, with $a,b,c,d>0$ (or a piece of the form $E_2^{(a)}E_1^{(b)}E_2^{(c)}E_1^{(d)}$, 
with $a,b,c,d>0$, which is done completely analogously). Then, we have that at least one of the inequalities $b<a+c$ or 
$c<b+d$ is satisfied. Suppose that the first one is satisfied (the second  one is done in the same way). Then by Proposition 
\ref{prgl} we have
$$E_1^{(a)}E_2^{(b)}E_1^{(c)}E_2^{(d)}=\sum_{{\scriptsize{\begin{array}{c} p+r=b \\ p\le c \\ r \le a\end{array}}}}
\left[\begin{array}{c}a+c-b\\ c-p \end{array}\right]  E_2^{(p)}E_1^{(a+c)}E_2^{(r)}E_2^{(d)}=$$
$$=\sum_{{\scriptsize{\begin{array}{c} p+r=b \\ p\le c \\ r \le a\end{array}}}}
\left[\begin{array}{c}a+c-b\\ c-p \end{array}\right]\left[\begin{array}{c}r+d\\ r \end{array}\right]  E_2^{(p)}E_1^{(a+c)}E_2^{(r+d)},$$
i.e. $v$ can be written as a linear combination of the monomials of length at most $l-1$, which proves the first part of the theorem.

As for the coefficients, they are all sums and products of the quantum binomial coefficients, and so they are from $\N[q,q^{-1}]$, as wanted. \kraj

\section{Category $\dot{\U}^+_n$}\label{sec3}

A categorification of the positive half of quantum $\mathfrak{sl}_n$ (and also of an arbitrary quantum group 
$U^+_q(\mathfrak{g})$) was defined in \cite{kl2}. The categorification there is given by means of the category of projective modules over certain (diagrammaticaly defined) rings $R(\nu)$ (see \cite{kl1}). In this paper we will use a categorification in terms of diagrammatic category $\mathcal{U}_n^+$, which is in the spirit of the categorification of the whole quantum $\mathfrak{sl}_n$ (see \cite{kl3}).\\

Before going to the definition of $\U_n^+$, first we recall some notation and explain the diagrams that appear in its definition (see also \cite{kl3}). 

Let $n\ge 2$ be fixed. By $\Seq$ we denote the set of all sequence of finite length of the numbers from $\{1,\ldots,n-1\}$. The entries of the elements of $\Seq$ we call \emph{colors}. For $\nu=(\nu_1,\ldots,\nu_{n-1})\in \N^{n-1}$, by $\Seq(\nu)$ we denote the set of all sequences $\underline{i}=(i_1,\ldots,i_k)\in \Seq$ such that $\sharp\{j| i_j=l\}=\nu_l$, for all $l=1,\ldots,n-1$.
Note that, in particular, we have that $k=\sum_l \nu_l$.\\

We will use the following diagrammatic calculus of planar diagrams: We consider 
collections of arcs on the plane connecting $k$ points on one horizontal line with $k$ 
points on another horizontal line. The positions of $k$ points on the horizontal line 
are always the points $\{1,\ldots,k\}\in\R$. Each arc is labelled by a number from the 
set $\{1,\ldots,n-1\}$ (called the \emph{color} of an arc). We require that arcs have no critical points when projected to $y$-axis. Also, arcs can intersect, but no triple intersections are allowed. Finally, an arc can carry dots.

The following is an example of a planar diagram:

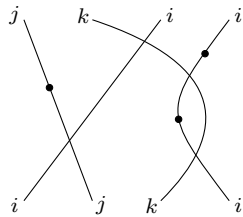
\begin{center}
\begin{tikzpicture} [scale=0.6]

\draw (0,-2)-- (3,2);
\draw (1.5,-2)--(0,2);

\draw (-0.2,-2.1) node {$\scriptstyle{i}$};
\draw (1.7,-2.1) node {$\scriptstyle{j}$};

\draw (-0.2,2.1) node {$\scriptstyle{j}$};
\draw (3.2,2.1) node {$\scriptstyle{i}$};

\draw (3,-2)..controls (4.5,-0.5) and (4.5,1) ..(1.5,2);
\draw (4.5,-2)..controls (3,0)..(4.5,2);

\draw (2.8, -2.1) node {$\scriptstyle{k}$};
\draw (1.3,2.1) node {$\scriptstyle{k}$};

\draw (4.7,-2.1) node {$\scriptstyle{i}$};
\draw (4.7,2.1) node {$\scriptstyle{i}$};

\filldraw[black] (0.56,0.5) circle (2pt)
(3.39,-0.2) circle (2pt)
(3.97,1.25) circle (2pt);

\end{tikzpicture}
\end{center}

We identify two planar diagrams if there exists a planar isotopy between them, that does 
not change the combinatorial type of the diagram and do not create critical points for 
the projection onto the $y$-axis.

Note that since we are not allowing the arcs to have critical points when projected to $y$-axis, we can assume that they are always oriented upwards. In particular, we can see a planar diagram as going from the sequence corresponding to the colors of the bottom end of the strands going to the sequence corresponding to the colors of the top end of the strands (we read the ends of the strands from left to right). For a diagram $D$, we denote the sequence that corresponds to the bottom (top) end of the strands by $\Bot(D)$ ($\Top(D)$, respectively). In the example above, we have $\Bot(D)=(i,j,k,i)$ and $\Top(D)=(j,k,i,i)$.

Each diagram has a degree, by setting that the degree of a dot is equal to $2$, while the degree of a crossing between two arcs that are colored by colors $i$ and $j$ is equal to $-i\cdot j$. In other words, for $i=j$ the degree of a crossing is equal to $-2$, for $|i-j|=1$ (adjacent colors) the degree of a crossing is equal to $1$, while for $|i-j|\ge 2$ (distant colors) the degree of a crossing is equal to $0$:
\begin{center}
\begin{tikzpicture} [scale=0.3]
\draw (0,-1.5)-- (0,1.5);

\draw (8,-1.5)--(10,1.5);
\draw (8,1.5)--(10,-1.5);

\draw (4,0) node {$,$};
\draw (-0.2,-1.7) node {$\scriptstyle{i}$};

\draw (7.8,-1.7) node {$\scriptstyle{i}$};
\draw (10.2,-1.7) node {$\scriptstyle{j}$};

\filldraw[black] (0,0) circle (4pt);

\draw (-6,-3.2) node {degree:};

\draw (0,-3.2)  node {$\scriptstyle{+2}$};

\draw (8.8,-3.2)  node {$\scriptstyle{-i\cdot j}$};

\end{tikzpicture}
\end{center}

We also use the following shorthand notation for the sequence of $d$ dots on an edge of a strand:
\begin{center}
\begin{tikzpicture} [scale=0.35]

\draw (0,-2)-- (0,2);
\draw (-0.2,-2.1) node {$\scriptstyle{i}$};
\draw (-0.4,0) node {$\scriptstyle{d}$};

\draw (1.2,0) node {$ := $};

\draw (3,-2) --(3,2);

\draw (2.8,-2.1) node {$\scriptstyle{i}$};
\draw (2.73,0.3) node {$\vdots$};

\draw (4,0) node {$\displaystyle \left. \right\} d $};

\filldraw[black] (0,0) circle (2pt)
(3,-1) circle (2pt)
(3,1) circle (2pt)
(3,-0.5) circle (2pt);

\end{tikzpicture}
\end{center}

\subsection{Category $\U^+_n$}
We define category $\U^+_n$ as follows:\\

$\U^+_n$ is the monoidal $\Z$-linear additive category whose objects and morphisms are the following:\\

$\bullet$    objects: for each $\underline{i}=(i_1,\ldots,i_k)\in \Seq$ and $t\in\Z$, we define $\E_{\underline{i}}\{t\}:=\E_{i_1}\ldots \E_{i_k}\{t\}$. An object of $\U^+_n$ is a formal finite direct sum of objects $\E_{\underline{i}}\{t\}$, with $i\in \Seq$ and $t\in\Z$.\\

$\bullet$    morphisms: for $\underline{i}=(i_1,\ldots,i_k)\in \Seq({\nu})$ and $\underline{j}=(j_1,\ldots,j_l)\in \Seq(\mu)$ the 
set $\Hom(\E_{\underline{i}}\{t\},\E_{\underline{j}}\{t'\})$ is empty, unless $\nu=\mu$. If $\nu=\mu$ (and consequently $k=l$), the 
morphisms from $\E_{\underline{i}}\{t\}$ to $\E_{\underline{j}}\{t'\}$ consists of finite $\Z$-linear combinations of planar diagrams going from $\underline{i}$ to $\underline{j}$, of degree $t-t'$, modulo the following set of local relations (note that all relations preserve the degree):

%Thus, the 2-morphisms are generated by the following generators:

%\begin{center}
%\begin{tikzpicture} [scale=0.3]
%\draw (0,-1.5)-- (0,1.5);

%\draw (8,-1.5)--(10,1.5);
%\draw (8,1.5)--(10,-1.5);

%\draw (4,0) node {$,$};
%\draw (-0.2,-1.7) node {$\scriptstyle{i}$};

%\draw (7.8,-1.7) node {$\scriptstyle{i}$};
%\draw (10.2,-1.7) node {$\scriptstyle{j}$};

%\filldraw[black] (0,0) circle (3pt);
%\end{tikzpicture}
%\end{center}

\begin{center}
\begin{tikzpicture} [scale=0.3]

\draw (6,-2)-- (8,2);
\draw (8,-2)--(6,2);

\draw (5.7,-2) node {$\scriptstyle{i}$};
\draw (8.3,-2) node {$\scriptstyle{i}$};

\draw (9,0) node {$\displaystyle{-}$};

\filldraw[black] (6.5,1) circle (3pt);

\draw (10,-2)-- (12,2);
\draw (12,-2)--(10,2);

\draw (9.7,-2) node {$\scriptstyle{i}$};
\draw (12.3,-2) node {$\scriptstyle{i}$};

\filldraw[black] (11.5,-1) circle (3pt);

\draw (13,0) node {$\displaystyle{=}$};

\draw (14,-2)-- (16,2);
\draw (16,-2)--(14,2);

\draw (13.7,-2) node {$\scriptstyle{i}$};
\draw (16.3,-2) node {$\scriptstyle{i}$};

\filldraw[black] (14.5,-1) circle (3pt);

\draw (17,0) node {$\displaystyle{-}$};

\draw (18,-2)-- (20,2);
\draw (20,-2)--(18,2);

\draw (17.7,-2) node {$\scriptstyle{i}$};
\draw (20.3,-2) node {$\scriptstyle{i}$};

\filldraw[black] (19.5,1) circle (3pt);

\draw (21,0) node {$\displaystyle{=}$};

\draw (22,-2)-- (22,2);
\draw (23.5,-2)--(23.5,2);

\draw (21.8,-2) node {$\scriptstyle{i}$};
\draw (23.7,-2) node {$\scriptstyle{i}$};

\end{tikzpicture}
\end{center}

\begin{center}
\begin{tikzpicture} [scale=0.3]

\draw (-8,-2)..controls (-6,0) ..(-8,2);
\draw (-6,-2)..controls (-8,0)..(-6,2);
\draw (-8.3,-2) node {$\scriptstyle{i}$};
\draw (-5.7,-2) node {$\scriptstyle{i}$};

\draw (-5,0) node {$\displaystyle{= 0} , $};

\draw (3.5,-2)--(6.5,2);
\draw (3.5,2)--(6.5,-2);
\draw (5,-2)..controls (3.5,0) ..(5,2);

\draw (3.7,-2) node {$\scriptstyle{i}$};
\draw (5.3,-2) node {$\scriptstyle{i}$};
\draw (6.8,-2) node {$\scriptstyle{i}$};

\draw (7.5,0) node {$=$};

\draw (8.5,-2)--(11.5,2);
\draw (8.5,2)--(11.5,-2);
\draw (10,-2)..controls (11.5,0) ..(10,2);

\draw (8.7,-2) node {$\scriptstyle{i}$};
\draw (10.3,-2) node {$\scriptstyle{i}$};
\draw (11.8,-2) node {$\scriptstyle{i}$};

\end{tikzpicture}
\end{center}

%\begin{center}
\begin{equation}\label{r2tan}
\begin{tikzpicture} [scale=0.3]
\draw (-1,-2)..controls (1,0) ..(-1,2);
\draw (1,-2)..controls (-1,0)..(1,2);
\draw (-1.3,-2) node {$\scriptstyle{i}$};
\draw (1.3,-2) node {$\scriptstyle{j}$};

\draw (2,0) node {$\displaystyle{=}$};

\draw (3,-2)--(3,2);
\draw (4.5,-2)--(4.5,2);

\draw (5.5,0) node {$+$};

\draw (6.5,-2)--(6.5,2);
\draw (8,-2)--(8,2);

\filldraw[black] (3,0) circle (3pt)
(8,0) circle (3pt);

\draw (15,0) node {,$\quad$ when $\quad i\cdot j=-1$};

\end{tikzpicture}
\end{equation}
%\end{center}

\begin{center}
\begin{tikzpicture}[scale=0.3]

\draw (11,-2)..controls (13,0) ..(11,2);
\draw (13,-2)..controls (11,0)..(13,2);
\draw (10.7,-2) node {$\scriptstyle{i}$};
\draw (13.3,-2) node {$\scriptstyle{j}$};

\draw (14,0) node {$\displaystyle{=}$};

\draw (15,-2)--(15,2);
\draw (16.5,-2)--(16.5,2);

\draw (14.7,-2) node {$\scriptstyle{i}$};
\draw (16.8,-2) node {$\scriptstyle{j}$};

\draw (23,0) node {,$\quad$ when $\quad i\cdot j=0$};

\end{tikzpicture}
\end{center}

\begin{center}
\begin{tikzpicture} [scale=0.3]
\draw (6,-2)-- (8,2);
\draw (8,-2)--(6,2);

\draw (5.7,-2) node {$\scriptstyle{i}$};
\draw (8.3,-2) node {$\scriptstyle{j}$};

\draw (9,0) node {$\displaystyle{=}$};

\filldraw[black] (6.5,1) circle (3pt);

\draw (10,-2)-- (12,2);
\draw (12,-2)--(10,2);

\draw (9.7,-2) node {$\scriptstyle{i}$};
\draw (12.3,-2) node {$\scriptstyle{j}$};

\filldraw[black] (11.5,-1) circle (3pt);

\draw (15,0) node {and};

%\draw (13,0) node {$\displaystyle{=}$};

\draw (19,-2)-- (21,2);
\draw (21,-2)--(19,2);

\draw (18.7,-2) node {$\scriptstyle{i}$};
\draw (21.3,-2) node {$\scriptstyle{j}$};

\filldraw[black] (19.5,-1) circle (3pt);

\draw (22,0) node {$\displaystyle{=}$};

\draw (23,-2)-- (25,2);
\draw (25,-2)--(23,2);

\draw (22.7,-2) node {$\scriptstyle{i}$};
\draw (25.3,-2) node {$\scriptstyle{j}$};

\filldraw[black] (24.5,1) circle (3pt);

\draw (31.5,0) node {, $\quad$ when $\quad i\ne j$};

\end{tikzpicture}
\end{center}

\begin{center}
\begin{tikzpicture} [scale=0.4]

\draw (-1.5,-2)--(1.5,2);
\draw (-1.5,2)--(1.5,-2);
\draw (0,-2)..controls (-1.5,0) ..(0,2);

\draw (-1.3,-2) node {$\scriptstyle{i}$};
\draw (0.3,-2) node {$\scriptstyle{j}$};
\draw (1.8,-2) node {$\scriptstyle{k}$};

\draw (2.5,0) node {$=$};

\draw (3.5,-2)--(6.5,2);
\draw (3.5,2)--(6.5,-2);
\draw (5,-2)..controls (6.5,0) ..(5,2);

\draw (3.7,-2) node {$\scriptstyle{i}$};
\draw (5.3,-2) node {$\scriptstyle{j}$};
\draw (6.8,-2) node {$\scriptstyle{k}$};

\draw (14,0) node {, if $i\ne k$ or $i\cdot j\ne -1$};

\end{tikzpicture}
\end{center}

%\begin{center}
\begin{equation}\label{r3tan}
\begin{tikzpicture} [scale=0.4]

\draw (-1.5,-2)--(1.5,2);
\draw (-1.5,2)--(1.5,-2);
\draw (0,-2)..controls (-1.5,0) ..(0,2);

\draw (-1.3,-2) node {$\scriptstyle{i}$};
\draw (0.3,-2) node {$\scriptstyle{j}$};
\draw (1.8,-2) node {$\scriptstyle{i}$};

\draw (2.5,0) node {$=$};

\draw (3.5,-2)--(6.5,2);
\draw (3.5,2)--(6.5,-2);
\draw (5,-2)..controls (6.5,0) ..(5,2);

\draw (3.7,-2) node {$\scriptstyle{i}$};
\draw (5.3,-2) node {$\scriptstyle{j}$};
\draw (6.8,-2) node {$\scriptstyle{i}$};

\draw (8,0) node {$+$};

\draw (9,-2)--(9,2);
\draw (10.5,-2)--(10.5,2);
\draw (12,-2)--(12,2);

\draw (9.3,-2) node {$\scriptstyle{i}$};
\draw (10.8,-2) node {$\scriptstyle{j}$};
\draw (12.3,-2) node {$\scriptstyle{i}$};

\draw (18,0) node {, if $i\cdot j= -1$};

\end{tikzpicture}
\end{equation}
%\end{center}

This ends the definition of $\U_n^+$. \\

We have the following relation in $\U_n^+$:

\begin{proposition}[Dot Migration]\cite{sl2}
We have

\begin{center}
\begin{tikzpicture} [scale=0.4]

\draw (6,-2)-- (8,2);
\draw (8,-2)--(6,2);

\draw (5.7,-2) node {$\scriptstyle{i}$};
\draw (8.3,-2) node {$\scriptstyle{i}$};

\draw (9,0) node {$\displaystyle{-}$};

\filldraw[black] (6.5,1) circle (3pt);
\draw (6.1,1) node {$\scriptstyle{d}$};

\draw (10,-2)-- (12,2);
\draw (12,-2)--(10,2);

\draw (9.7,-2) node {$\scriptstyle{i}$};
\draw (12.3,-2) node {$\scriptstyle{i}$};

\filldraw[black] (11.5,-1) circle (3pt);
\draw (11.9,-1) node {$\scriptstyle{d}$};

\draw (13,0) node {$\displaystyle{=}$};

\draw (14,-2)-- (16,2);
\draw (16,-2)--(14,2);

\draw (13.7,-2) node {$\scriptstyle{i}$};
\draw (16.3,-2) node {$\scriptstyle{i}$};

\filldraw[black] (14.5,-1) circle (3pt);
\draw (14.1,-1) node {$\scriptstyle{d}$};

\draw (17,0) node {$\displaystyle{-}$};

\draw (18,-2)-- (20,2);
\draw (20,-2)--(18,2);

\draw (17.7,-2) node {$\scriptstyle{i}$};
\draw (20.3,-2) node {$\scriptstyle{i}$};

\filldraw[black] (19.5,1) circle (3pt);
\draw (19.9,1) node {$\scriptstyle{d}$};

\draw (22.5,0) node {$\displaystyle{=} \sum_{r+s=d-1}$};

\draw (25.5,-2)-- (25.5,2);
\draw (27,-2)--(27,2);

\draw (25.3,-2) node {$\scriptstyle{i}$};
\draw (27.2,-2) node {$\scriptstyle{i}$};

\filldraw[black] (25.5,0) circle (3pt);
\draw (25.2,0) node {$\scriptstyle{r}$};

\filldraw[black] (27,0) circle (3pt);
\draw (27.4,0) node {$\scriptstyle{s}$};

\end{tikzpicture}
\end{center}

\end{proposition}

\subsection{Category $\dot{\U}_n^+$ and thick calculus}

In \cite{thick}, the extension of the calculus to thick edges have been introduced. Thick lines categorify the divided powers $E_i^{(a)}$ (see below and Section 4 of \cite{thick}). 

A thick line is defined in terms of ``ordinary" lines from above, and drawn as a strand with an additional label (natural number) $a$ (also called the \emph{thickness} of a strand). In particular, the ordinary strands from above correspond to the case $a=1$, and are also called \emph{thin} edges or thin strands.  
For any color the thick edges are defined in the same way as in \cite{thick}, and we refer the reader to that paper for more details. Here we just recall the basic facts that will be used later on.

Also, now the trivalent vertices in planar diagrams are allowed (called Splitters in \cite{thick}), such the sum of the thicknesses of the incoming edges is equal to the sum of the thicknesses of the outgoing edges. Recall that we are assuming that all strands are oriented upwards. The degree of a trivalent vertex (for any color - the labels on the pictures below represent thicknesses)

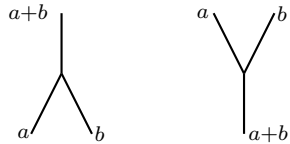
\begin{center}
\begin{tikzpicture} [scale=0.4]

\draw[thick] (0,-2)-- (1,0);
\draw[thick] (2,-2)--(1,0);
\draw[thick] (1,0)--(1,2);

\draw (-0.25,-2) node {$\scriptstyle{a}$};
\draw (2.25,-2) node {$\scriptstyle{b}$};
\draw (-0.1,2) node {$\scriptstyle{a+b}$};

\draw[thick] (7,-2)-- (7,0);
\draw[thick] (6,2)--(7,0);
\draw[thick] (7,0)--(8,2);

\draw (7.8,-2) node {$\scriptstyle{a+b}$};
\draw (5.65,2) node {$\scriptstyle{a}$};
\draw (8.25,2) node {$\scriptstyle{b}$};
\end{tikzpicture}
\end{center}

\noindent is equal to $-ab$.\\

\subsubsection{Category $\dot{\U}_n^+$}

For a category $\mathcal{C}$, the Karoubi envelope $Kar(\mathcal{C})$ is the smallest category containing $\mathcal{C}$, such that all idempotents split (for more details, see e.g. Section 3.4 of \cite{thick}).

We define the category $\dot{\U}^+_n$ as the Karoubi envelope of the category $\U^+_n$. \\

Now, in the category $\dot{\U}^+_n$, the planar diagrams with thick edges from above can be interpreted as morphisms whose bottom and top end correspond to certain objects of $\dot{\U}^+_n$. In particular, the object corresponding to bottom (or the top end)  of an arc of color $i$ and thickness $a$ is denoted $\E_i^{(a)}$. 

As in \cite{thick} (see also \cite{thickn}), the category $\dot{\U}^+_n$ categorifies $U_q^+(\mathfrak{sl}_n)$, in a sense that its split Grothendieck group is isomorphic to the integral form of  $U_q^+(\mathfrak{sl}_n)$. The isomorphism sends the class of  $\E_i^{(a)}$ to the generator $E_i^{(a)}$ of $U_q^+(\mathfrak{sl}_n)$.

\vskip 0.5cm

\subsubsection{Some properties of thick calculus}

Below we give some of the basic properties of thick edges that we shall use in this paper (see \cite{thick} for more details).
Note that the labels of the strands below denote thickness.

\begin{proposition}[Associativity of splitters]
For arbitrary color $i$ (drawn as thick line), we have the following:

\begin{center}
\begin{tikzpicture}[scale=0.4]
\draw [thick] (-2,2)--(0,-0.5);
\draw [thick](-1,0.75)--(0,2);
\draw [thick](0,-0.5)--(2,2);
\draw [thick](0,-0.5)--(0,-2);

\draw (-2.3,2) node {$\scriptstyle{a}$};
\draw (-0.45,2) node {$\scriptstyle{b}$};
\draw (2.3,2) node {$\scriptstyle{c}$};
\draw (-1.3,0) node {$\scriptstyle{a+b}$};
\draw (1.1,-2.2) node {$\scriptstyle{a+b+c}$};

\draw (3,0) node {$=$};

\draw [thick](4,2)--(6,-0.5);
\draw [thick](7,0.75)--(6,2);
\draw [thick](6,-0.5)--(8,2);
\draw [thick](6,-0.5)--(6,-2);

\draw (3.7,2) node {$\scriptstyle{a}$};
\draw (5.7,2) node {$\scriptstyle{b}$};
\draw (8.3,2) node {$\scriptstyle{c}$};
\draw (7.4,0) node {$\scriptstyle{b+c}$};
\draw (7.5,-2.2) node {$\scriptstyle{a+b+c}$};

\end{tikzpicture}
\end{center}

\end{proposition}

\begin{proposition}[Pitchfork lemma]
For any two colors $i$ (drawn as thick line) and $j$ (drawn dashed), we have:
\begin{center}
\begin{tikzpicture}[scale=0.5]
\draw [thick](-2,-2)--(-1,-1);
\draw [thick](-1,-1)..controls (-0.5,0.5)..(1,2);
\draw [thick](-1,-1)..controls (0.5,-0.5)..(2,1);
\draw [thick,dashed] (-2,2)--(2,-2);

\draw (1.3,2) node {$\scriptstyle a$};
\draw (-2,-2.35) node {$\scriptstyle a+b$};
\draw (2.4,1.2) node {$\scriptstyle b$};
\draw (2.3,-2) node {$\scriptstyle c$};

\draw (3,0) node {$=$};

\draw [thick](4,-2)--(6.5,0.5);
\draw [thick](6.5,0.5)..controls (6.6,1.5)..(7,2);
\draw [thick](6.5,0.5)..controls (7.5,0.6)..(8,1);

\draw [thick,dashed] (4,2)--(8,-2);

\draw (7.4,2) node {$\scriptstyle a$};
\draw (4.5,-2.3) node {$\scriptstyle a+b$};
\draw (8.35,1.2) node {$\scriptstyle b$};
\draw (8.3,-2) node {$\scriptstyle c$};

\end{tikzpicture}
\end{center}

\end{proposition}

\begin{proposition}[Opening of a Thick Edge]
For any color $i$ (drawn as thick line), we have:

\begin{center}
\begin{tikzpicture}[scale=0.9]

\draw [thick](0,-1.5)--(1,-0.6);
\draw [thick](1,-0.6)--(1,0.6);
\draw [thick](1,0.6)--(0,1.5);
\draw [thick](2,-1.5)--(1,-0.6);
\draw [thick](1,0.6)--(2,1.5);
\draw (3,0) node {$=$};

\draw [thick](4,-1.5)--(7,1.5);
\draw [thick](4,1.5)--(7,-1.5);
\draw [thick](4.6,-0.91)--(4.6,0.91);

\draw (-0.5,1.6) node {$\scriptstyle{b+x}$};
\draw (-0.5,-1.6) node {$\scriptstyle{a+x}$};
\draw (0.3,0) node {$\scriptstyle{a+b+x}$};
\draw (2.2,1.6) node {$\scriptstyle{a}$};
\draw (2.2,-1.6) node {$\scriptstyle{b}$};

\draw (3.5,1.6) node {$\scriptstyle{b+x}$};
\draw (3.5,-1.6) node {$\scriptstyle{a+x}$};
\draw (4.3,0) node {$\scriptstyle{x}$};
\draw (7.2,1.6) node {$\scriptstyle{a}$};
\draw (7.2,-1.55) node {$\scriptstyle{b}$};
\draw (5,0.7) node {$\scriptstyle{b}$};
\draw (5,-0.7) node {$\scriptstyle{a}$};

\end{tikzpicture}
\end{center}

\end{proposition}

\subsection{Schur polynomials and labels on thick lines}

Thick lines are labelled with symmetric polynomials, that correspond to symmetric polynomials in dots on thin edges involved in the definition of a thick line (for precise definition see \cite{thick}). 
In particular, we label thick edges by Schur polynomials, which form the additive basis of the ring of symmetric polynomials. 

\subsubsection{Schur polynomials}
Here we recall briefly the definition and some basic notation and properties of the Schur polynomials. For more details on them see e.g. \cite{thick} or \cite{fulton}.

By a partition $\alpha=(\alpha_1,\ldots,\alpha_k)$, we mean a non-increasing sequence of non-negative integers. We identify two partitions if they differ by a sequence of zeros at the end. We set $|\alpha|=\sum_i \alpha_i$. If for some $a$ we have $\alpha_{a+1}=0$, we say that $\alpha$ has at most $a$ parts. We denote the set of all partitions with at most $a$ parts by 
$P(a)$. Furthermore, by $P(a,b)$ we denote the subset of all partitions $\alpha$ from $P(a)$ such that $\alpha_1\le b$.

By $\bar{\alpha}$ we denote the dual (conjugate) partition of $\alpha$, i.e. $\alpha_j=\sharp\{i| \alpha_i\ge j\}$. If $\alpha\in P(a,b)$, we define partition $\hat{\alpha}$ by $\hat{\alpha}=\overline{(b-\alpha_a,\ldots,b-\alpha_1)}$. Note that if $\alpha\in P(a,b)$, then $\bar{\alpha}\in P(b,a)$ and $\hat{\alpha}\in P(b,a)$.

For any partition $\alpha\in P(a)$, the Schur polynomial $\pi_{\alpha}$ is given by the formula:
\[
\pi_{\alpha}(x_1,x_2,\ldots,x_a)=\frac{|x_i^{\alpha_j+a-j}|}{\Delta},
\]
where $\Delta=\prod_{1\le r \le s\le a} (x_r-x_s)$, and $|x_i^{\alpha_j+a-j}|$ is the determinant of the $a\times a$ matrix whose $(i,j)$ entry is $x_i^{\alpha_j+a-j}$. We let $\pi_{\alpha}(x_1,x_2,\ldots,x_a)=0$ is some entry of $\alpha$ is negative ($\alpha$ is not a partition then), or if $\alpha_{a+1}>0$.

For three partitions $\alpha$, $\beta$ and $\gamma$, the Littlewood-Richardson coefficient $c_{\alpha,\beta}^{\gamma}$ are given by:
\[
\pi_{\alpha}\pi_{\beta}=\sum_{\gamma} c_{\alpha,\beta}^{\gamma} \pi_{\gamma}.
\]
The coefficients $c_{\alpha,\beta}^{\gamma}$ are nonnegative integers that can be nonzero only when $|\gamma|=|\alpha|+|\beta|$.

The Littlewood-Richardson coefficients can be naturally extended for more than three partitions: for partitions $\alpha_1,\ldots,\alpha_k$ and $\beta$, with $k\ge 2$, we define $c_{\alpha_1,\ldots,\alpha_k}^{\beta}$ by:
\[
\pi_{\alpha_1}\ldots\pi_{\alpha_k}=\sum_{\beta} c_{\alpha_1,\ldots,\alpha_k}^{\beta} \pi_{\beta}.
\]
For two partitions $\alpha$ and $\gamma$, we say that $\alpha\subset\gamma$, if $\alpha_i\le \gamma_i$ for all $i\ge 1$. 
If $\alpha\subset \gamma$, we define skew-Schur polynomial $\pi_{^{\gamma}/_{\alpha}}$ by:
\[
\pi_{^{\gamma}/_{\alpha}}=\sum_{\beta} c_{\alpha,\beta}^{\gamma} \pi_{\beta}.
\]

Finally, the elementary symmetric polynomials $\varepsilon_m(x_1,\ldots,x_a)$, for $m=0,\ldots,a$,  are special Schur polynomials: $\varepsilon_m(x_1,\ldots,x_a)=\pi_{(\underbrace{1,1,\ldots,1}_m)} ((x_1,\ldots,x_a)$.

\subsubsection{Thick edges with labels}
Going back to the labels of thick edges, the thick line of thickness $a$ can be labelled with the Schur polynomials $\pi_{\alpha}$, with $\alpha\in P(a)$. If $\alpha\notin P(a)$, the diagram is equal to zero.

In terms of thin edges and dots, Schur polynomial $\pi_{\alpha}$, with $\alpha=(\alpha_1,\alpha_2,\ldots,\alpha_a)$, labelling thick edge of thickness $a$ is given by:
\begin{center}
\begin{tikzpicture}[scale=0.7]

\draw [very thick] (0,1.5)--(0,-1.5);

\draw (-0.5,0) node {$\scriptstyle{\pi_{\alpha}}$};
\draw (2.1,0.15) node {$\scriptscriptstyle{\alpha_{1}+a-1}$};
\draw (4.25,0.2) node {$\scriptscriptstyle{\alpha_{2}\!+\!a\!-\!2}$};
\draw (5.65,-0.3) node {$\scriptscriptstyle{\alpha_{a-1}\!+\!1}$};
\draw (7.7,0) node {$\scriptstyle{\alpha_{a}}$};
\draw (5,0) node {$\cdots$};

%\draw [snake=snake,segment amplitude=.2mm,segment length=2mm](-15,-1.5)--(-15,1.5);

\filldraw[black] (0,0) circle (1.8pt)
                  
                 (6.5,0) circle (1.8pt)
                 (7.25,0) circle (1.8pt)
                 (3.522,0) circle (1.8pt)
                 (2.725,0) circle (1.8pt);

\draw (0.2,-1.5) node {$\scriptstyle{a}$};

\draw (5.2,2) node {$\scriptstyle{a}$};
\draw (5.2,-2) node {$\scriptstyle{a}$};

\draw (1,0) node{$=$};
\draw [very thick] (5,-2)--(5,-1.5);
\draw [very thick] (5,2)--(5,1.5);
\draw (5,1.5)..controls (2,0.75) and (2,-0.75)..(5,-1.5);
\draw (5,1.5)..controls (8,0.75) and (8,-0.75)..(5,-1.5);
\draw (5,1.5)..controls (3,0.75) and (3,-0.75)..(5,-1.5);
\draw (5,1.5)..controls (7,0.75) and (7,-0.75)..(5,-1.5);
\end{tikzpicture}
\end{center}

By ``exploding" a thick edge into thin edges, we obtain the diagrams that are antisymmetric with respect to the exchange of dots on two neighboring strands:
\begin{center}
\begin{tikzpicture}[scale=0.5]

\draw [very thick] (11,-2)--(11,-1.2);
\draw [very thick] (11,2)--(11,1.2);
\draw (11,1.2)..controls (10,0)..(11,-1.2);
\draw (11,1.2)..controls (12,0)..(11,-1.2);
\draw (10,0) node {$\scriptstyle{a}$};
\draw (12,0) node {$\scriptstyle{b}$};
\draw (13,0) node {$= - $};
\draw [very thick] (15,-2)--(15,-1.2);
\draw [very thick] (15,2)--(15,1.2);
\draw (15,1.2)..controls (14,0)..(15,-1.2);
\draw (15,1.2)..controls (16,0)..(15,-1.2);
\draw (14,0) node {$\scriptstyle{b}$};
\draw (16,0) node {$\scriptstyle{a}$};

\filldraw[black] (10.23,0) circle (1.7pt)
                 (11.75,0) circle (1.7pt)
                 (14.23,0) circle (1.7pt)
                 (15.75,0) circle (1.7pt);
                 %(3.45,0) circle (1.5pt);

\end{tikzpicture}
\end{center}

Also, we have the following:

\begin{lemma} {}$\quad$

\begin{center}

\begin{tikzpicture}[scale=0.7]
%\draw (-0.5,0) node {$\scriptstyle{\pi_{\alpha}}$};
\draw (2.1,0) node {$\scriptscriptstyle{x_1}$};
\draw (4,0) node {$\scriptscriptstyle{x_2}$};
\draw (5.8,0) node {$\scriptscriptstyle{x_{a-1}}$};
\draw (7.6,0) node {$\scriptstyle{x_{a}}$};
\draw (5,0) node {$\cdots$};

%\draw [snake=snake,segment amplitude=.2mm,segment length=2mm](-15,-1.5)--(-15,1.5);

\filldraw[black] %(0,0) circle (1.8pt)
                  (7.23,0) circle (1.8pt)
                 (6.5,0) circle (1.8pt)
                 
(3.522,0) circle (1.8pt)
                 (2.725,0) circle (1.8pt);

%\draw (0.2,-1.5) node {$\scriptstyle{a}$};

\draw (5.2,2) node {$\scriptstyle{a}$};
\draw (5.2,-2) node {$\scriptstyle{a}$};

%\draw (1,0) node{$=$};
\draw [very thick] (5,-2)--(5,-1.5);
\draw [very thick] (5,2)--(5,1.5);
\draw (5,1.5)..controls (2,0.75) and (2,-0.75)..(5,-1.5);
\draw (5,1.5)..controls (8,0.75) and (8,-0.75)..(5,-1.5);
\draw (5,1.5)..controls (3,0.75) and (3,-0.75)..(5,-1.5);
\draw (5,1.5)..controls (7,0.75) and (7,-0.75)..(5,-1.5);

\draw (12.3,0) node {$\displaystyle{\ne 0}\quad\Rightarrow\quad x_r\ne x_s , \textrm{ for all } r\ne s$};
\end{tikzpicture}
\end{center}

Moreover, if $\max_i\{x_i\}=a-1$, then the diagram from above can be nonzero if and only if there exists a permutation $\sigma$ of $\{0,1,\ldots,a-1\}$, such that $x_i=\sigma_i$, $i=0,\ldots,a-1$, and then  

\begin{center}
\begin{tikzpicture}[scale=0.55]
%\draw (-0.5,0) node {$\scriptstyle{\pi_{\alpha}}$};
\draw (2.1,0) node {$\scriptscriptstyle{x_1}$};
\draw (4,0) node {$\scriptscriptstyle{x_2}$};
\draw (5.8,0) node {$\scriptscriptstyle{x_{a-1}}$};
\draw (7.6,0) node {$\scriptstyle{x_{a}}$};
\draw (5,0) node {$\cdots$};
%\draw [snake=snake,segment amplitude=.2mm,segment length=2mm](-15,-1.5)--(-15,1.5);
\filldraw[black] %(0,0) circle (1.8pt)
                  (7.23,0) circle (1.8pt)
                 (6.5,0) circle (1.8pt)
                 
(3.522,0) circle (1.8pt)
                 (2.725,0) circle (1.8pt);

%\draw (0.2,-1.5) node {$\scriptstyle{a}$};

\draw (5.2,2) node {$\scriptstyle{a}$};
\draw (5.2,-2) node {$\scriptstyle{a}$};

%\draw (1,0) node{$=$};
\draw [very thick] (5,-2)--(5,-1.5);
\draw [very thick] (5,2)--(5,1.5);
\draw (5,1.5)..controls (2,0.75) and (2,-0.75)..(5,-1.5);
\draw (5,1.5)..controls (8,0.75) and (8,-0.75)..(5,-1.5);
\draw (5,1.5)..controls (3,0.75) and (3,-0.75)..(5,-1.5);
\draw (5,1.5)..controls (7,0.75) and (7,-0.75)..(5,-1.5);

\draw (10,0) node {$\displaystyle{=\sgn \sigma}$};

\draw [very thick] (12,-2)--(12,2);
\draw (12.2,-2.1) node {$\scriptstyle{a}$}; 
\end{tikzpicture}
\end{center}
\end{lemma}

The above lemma implies the following

\begin{lemma}\label{vazlem}
Let $\alpha\in P(a,x)$ and $\beta\in P(b,y)$ be partitions. Then we have that 
\begin{center}
\begin{tikzpicture}[scale=0.5]
\draw [very thick] (11,-2)--(11,-1.2);
\draw [very thick] (11,2)--(11,1.2);
\draw (11,1.2)..controls (10,0)..(11,-1.2);
\draw (11,1.2)..controls (12,0)..(11,-1.2);
\draw (9.7,0) node {$\scriptstyle{\pi_{\alpha}}$};
\draw (12.3,0) node {$\scriptstyle{\pi_{\beta}}$};
\draw (13.5,0) node {$= s$};
\draw [very thick] (15,-2)--(15,2);

\draw (15.4,-2.2) node {$\scriptstyle{a+b}$};
\draw (15.55,0) node {$\scriptstyle{\pi_{\gamma}}$};
%\draw (15,2)--(15,1.2);
%\draw (15,1.2)..controls (14,0)..(15,-1.2);
%\draw (15,1.2)..controls (16,0)..(15,-1.2);
%\draw (14,0) node {$b$};
%\draw (15.5,0) node {$a+b$};
\draw (10.6,-2.3) node {$\scriptstyle{a+b}$};
\draw (10.6,2.3) node {$\scriptstyle{a+b}$};

\draw (10.5,-1) node {$\scriptstyle{a}$};
\draw (11.5,-1) node {$\scriptstyle{b}$};

\filldraw[black] (15,0) circle (2pt);
\filldraw[black] (10.2,0) circle (2pt);
\filldraw[black] (11.8,0) circle (2pt);

\end{tikzpicture}
\end{center}
\noindent for some partition $\gamma\in P(a+b,x-b+y-a)$ and $s\in\{-1,0,1\}$. If $s\ne 0$, then $|\gamma|=|\alpha|+|\beta|-ab$.

Moreover, if $\alpha\in P(a,b)$ and $\beta\in P(b,a)$, then we have:
\begin{center}
\begin{tikzpicture}[scale=0.5]
\draw [very thick] (9.5,-2)--(9.5,-1.2);
\draw [very thick] (9.5,2)--(9.5,1.2);
\draw (9.5,1.2)..controls (8.5,0)..(9.5,-1.2);
\draw (9.5,1.2)..controls (10.5,0)..(9.5,-1.2);
\draw (8.2,0) node {$\scriptstyle{\pi_{\alpha}}$};
\draw (10.8,0) node {$\scriptstyle{\pi_{\beta}}$};
\draw (13.5,0) node {$= \delta_{\beta,\hat{\alpha}} (-1)^{|\beta|}$};
\draw [very thick] (16.2,-2)--(16.2,2);

\draw (16.4,-2.2) node {$\scriptstyle{a+b}$};
%\draw (15.4,0) node {$\scriptstyle{\pi_{\gamma}}$};
%\draw (15,2)--(15,1.2);
%\draw (15,1.2)..controls (14,0)..(15,-1.2);
%\draw (15,1.2)..controls (16,0)..(15,-1.2);
%\draw (14,0) node {$b$};
%\draw (15.5,0) node {$a+b$};
\draw (10,-2.3) node {$\scriptstyle{a+b}$};
\draw (10,2.3) node {$\scriptstyle{a+b}$};

\draw (9,-1) node {$\scriptstyle{a}$};
\draw (10,-1) node {$\scriptstyle{b}$};

%\filldraw[black] (15,0) circle (2pt);
\filldraw[black] (8.7,0) circle (2pt);
\filldraw[black] (10.3,0) circle (2pt);

\end{tikzpicture}
\end{center}
In particular, the left hand side can be nonzero only when $|\alpha|+|\beta|=ab$.
%then the right-hand side is nonzero only if $\beta=\hat{\alpha}$ and in that case $\gamma=0$ and $s=(-1)^{|\alpha|}$. 
\end{lemma}

\begin{lemma}\label{pomocna}
Let $\gamma\in P(a)$ and $\psi\in P(a,b)$ be partitions, and let $K=(\underbrace{b,\ldots,b}_{a})$. Then
\begin{center}
\begin{tikzpicture}[scale=0.5]
\draw [very thick] (11,-2)--(11,-1.2);
\draw [very thick] (11,2)--(11,1.2);
\draw (11,1.2)..controls (10,0)..(11,-1.2);
\draw (11,1.2)..controls (12,0)..(11,-1.2);
\draw (9.7,-0.4) node {$\scriptstyle{\pi_{\gamma}}$};
\draw (9.7,0.4) node {$\scriptstyle{\pi_{\psi}}$};
\draw (13.7,0) node {$=$};
\draw [very thick] (15,-2)--(15,2);

\draw (15.4,-2.2) node {$\scriptstyle{a+b}$};
\draw (16.5,0) node {$\scriptstyle{\pi_{_{^{\gamma}/_{(K-\psi)}}}}$};
%\draw (15,2)--(15,1.2);
%\draw (15,1.2)..controls (14,0)..(15,-1.2);
%\draw (15,1.2)..controls (16,0)..(15,-1.2);
%\draw (14,0) node {$b$};
%\draw (15.5,0) node {$a+b$};
\draw (10.6,-2.3) node {$\scriptstyle{a+b}$};
\draw (10.6,2.3) node {$\scriptstyle{a+b}$};

\draw (10.5,-1) node {$\scriptstyle{a}$};
\draw (11.5,-1) node {$\scriptstyle{b}$};

\filldraw[black] (15,0) circle (2pt);
\filldraw[black] (10.33,-0.3) circle (2pt);
\filldraw[black] (10.33,0.3) circle (2pt);

%\filldraw[black] (11.8,0) circle (2pt);

\end{tikzpicture}
\end{center}
where $K-\psi=(b-\psi_a,b-\psi_{a-1},\ldots,b-\psi_1)$.
\end{lemma}

\textbf{Proof:}\\

\begin{center}
\begin{tikzpicture}[scale=0.5]
\draw [very thick] (11,-2)--(11,-1.2);
\draw [very thick] (11,2)--(11,1.2);
\draw (11,1.2)..controls (10,0)..(11,-1.2);
\draw (11,1.2)..controls (12,0)..(11,-1.2);
\draw (9.7,-0.4) node {$\scriptstyle{\pi_{\gamma}}$};
\draw (9.7,0.4) node {$\scriptstyle{\pi_{\psi}}$};
\draw (14.8,0) node {$\displaystyle{=\sum_{\mu\in P(a)} c^{\mu}_{\gamma,\psi}}$};

\draw [very thick] (19,-2)--(19,-1.2);
\draw [very thick] (19,2)--(19,1.2);
\draw (19,1.2)..controls (18,0)..(19,-1.2);
\draw (19,1.2)..controls (20,0)..(19,-1.2);
\draw (17.7,0.4) node {$\scriptstyle{\pi_{\mu}}$};

\draw (18.5,-1) node {$\scriptstyle{a}$};
\draw (19.5,-1) node {$\scriptstyle{b}$};
\draw (18.6,2.3) node {$\scriptstyle{a+b}$};

%\draw (15.4,-2.2) node {$\scriptstyle{a+b}$};
%\draw (16.5,0) node {$\scriptstyle{\pi_{^{\gamma}/_{(K-\psi)}}}$};
%\draw (15,2)--(15,1.2);
%\draw (15,1.2)..controls (14,0)..(15,-1.2);
%\draw (15,1.2)..controls (16,0)..(15,-1.2);
%\draw (14,0) node {$b$};
%\draw (15.5,0) node {$a+b$};
\draw (10.6,-2.3) node {$\scriptstyle{a+b}$};
\draw (10.6,2.3) node {$\scriptstyle{a+b}$};

\draw (10.5,-1) node {$\scriptstyle{a}$};
\draw (11.5,-1) node {$\scriptstyle{b}$};

%\filldraw[black] (15,0) circle (2pt);
\filldraw[black] (10.33,-0.3) circle (2pt);
\filldraw[black] (10.33,0.3) circle (2pt);
\filldraw[black] (18.33,0.3) circle (2pt);
\filldraw[black] (26,0) circle (2pt);
 \draw (23,0) node {$\displaystyle{=\sum_{\nu\in P(a)} c^{\nu+K}_{\gamma,\psi}}$}; 
\draw [very thick] (26,-2)--(26,2);
\draw (26.2,-2.1) node {$\scriptstyle{a+b}$};                
\draw (26.5,0) node {$\scriptstyle{\pi_\nu}$};

\end{tikzpicture}
\end{center}
\noindent where $\nu+K=(\nu_1+b,\ldots,\nu_a+b)$. Since $\pi_{K/\gamma}=\pi_{K-\gamma}$, we have $c^{\nu+k}_{\gamma,\psi}=c^{\gamma}_{\nu,K-\psi}$, and so $\pi_{\gamma/{(K-\psi)}}=\sum_{\nu} c_{\nu,K-\psi}^{\gamma} \pi_{\nu} = 
\sum_{\nu} c_{\gamma,\psi}^{\nu+K} \pi_{\nu}$, which gives the above lemma. \kraj\\

We shall also need the extensions of R2 and R3 like moves for thick edges, in the case when colors $i$ and $j$ of strands involved satisfy $i \cdot j=-1$.\\

Since in thick calculus each line carries two indices, we shall use the following convention for drawing diagrams, in order to have as small number labels on diagrams as possible.\\

{\bf {Notation convention:}} From now on, for two colors (indexes) that satisfy $i\cdot j=-1$, we shall draw lines corresponding to $i$ as normal (straight) lines, while the lines corresponding to color $j$, we shall draw as curly lines:
\begin{center}
\begin{tikzpicture}[scale=0.55]

\draw (0,0) node {$\displaystyle{\mathcal{E}_i^{(a)}}:$};
\draw [dashed, thin] (1,0)--(3,0);

\draw [thick] (2,-0.7)--(2,0.7);

\draw (1.8,-0.8) node {$\scriptstyle{a}$};

\draw (8,0) node {$\displaystyle{\mathcal{E}_j^{(b)}}:$};
\draw [dashed,thin] (9,0)--(11,0);

\draw [snake=snake,segment amplitude=.2mm,segment length=2mm,thick]
(10,-0.7)--(10,0.7);
 
 \draw (9.8,-0.8) node {$\scriptstyle{b}$};
\end{tikzpicture}
\end{center}

Thus, from now on, each line carries one label, and that label represents its thickness.\\

The two proposition below are straightforward extensions of the thin R2 and R3 relations (\ref{r2tan}) and (\ref{r3tan}) (see \cite{thick} and \cite{thickn} for more details).

\begin{proposition}[Thick R2 Move]
We have

\begin{center}
\begin{tikzpicture} [scale=0.5]
\draw [thick] (-1,-3)--(1,0);
\draw [thick] (1,0)--(-1,3);
\draw [thick,snake=snake,segment amplitude=.2mm,segment length=2mm](1,-3)--(-1,0);
\draw [thick,snake=snake,segment amplitude=.2mm,segment length=2mm](-1,0)--(1,3);
\draw (-1.3,-3) node {$\scriptstyle{a}$};
\draw (1.3,-3) node {$\scriptstyle{b}$};
\draw [thick] (3.7,0) node {$= \displaystyle{\sum_{\alpha\in P(a,b)}}$};
\draw [thick](6.5,-3)--(6.5,3);
\draw [thick,snake=snake,segment amplitude=.2mm,segment length=2mm](8,-3)--(8,3);
\draw (5.9,0) node {$\scriptstyle{\pi_{\alpha}}$};
\draw (8.55,0) node {$\scriptstyle{\pi_{\hat{\alpha}}}$};
\draw (6.3,-3.2) node {$\scriptstyle{a}$};
\draw (8.1,-3.2) node {$\scriptstyle{b}$};
\filldraw[black] (6.5,0) circle (2pt) 
(8,0) circle (2pt);
\end{tikzpicture}
\end{center}
\end{proposition}

\begin{proposition}[Thick R3 Move]
We have

\begin{center}
\begin{tikzpicture} [scale=0.5]
\draw [thick] (-1,-3)--(2,3);
\draw [thick](2,-3)--(-1,3);
\draw [thick, snake=snake,segment amplitude=.2mm,segment length=2mm](0.5,-3)--(-0.5,0);
\draw [thick, snake=snake,segment amplitude=.2mm,segment length=2mm](-0.5,0)--(0.5,3);

\draw (-1.3,-3) node {$\scriptstyle{a}$};
\draw (0.8,-3) node {$\scriptstyle{c}$};
\draw (2.2,-3) node {$\scriptstyle{b}$};

\draw (7,0) node {$= \displaystyle{\sum_{i=0}^{\min(a,b,c)}\sum_{\alpha,\beta,\gamma\in P(i,c-i)}}c^{K_i}_{\alpha\beta\gamma}$};
\draw [thick](12,-3)--(13,-1.5);
\draw [thick](12,3)--(13,1.5);
\draw [thick](18,-3)--(17,-1.5);
\draw [thick](18,3)--(17,1.5);

\draw [thick](13,1.5)--(17,-1.5);
\draw [thick](13,-1.5)--(17,1.5);
\draw [thick](13,1.5)--(13,-1.5);
\draw [thick](17,1.5)--(17,-1.5);

\draw [thick,snake=snake,segment amplitude=.2mm,segment length=2mm](15,-3)--(16,0);
\draw [thick,snake=snake,segment amplitude=.2mm,segment length=2mm](16,0)--(15,3);

\draw (12.4,3) node {$\scriptstyle{b}$};
\draw (18.2,3) node {$\scriptstyle{a}$};
\draw (11.7,-3) node {$\scriptstyle{a}$};
\draw (18.2,-3) node {$\scriptstyle{b}$};
\draw (15.3,-3) node {$\scriptstyle{c}$};

\draw (12.8,1.2) node {$\scriptstyle{i}$};
\draw (17.2,1.2) node {$\scriptstyle{i}$};

\draw (14.2,1.2) node {$\scriptstyle{b-i}$};
\draw (14,-1.2) node {$\scriptstyle{a-i}$};

\draw (12.4,0) node {$\scriptstyle{\pi_{\alpha}}$};
\draw (17.6,0) node {$\scriptstyle{\pi_{\beta}}$};
\draw (16,-2) node {$\scriptstyle{\pi_{\bar{\gamma}}}$};

\filldraw[black] (13,0) circle (2pt) 
(17,0) circle (2pt)
(15.35,-2) circle (2pt);
\end{tikzpicture}
\end{center}
where $K_i=(\underbrace{c-i,c-i,\ldots,c-i}_i)$, for $i>0$, and $K_0=0$.
\end{proposition}

\section{Indecomposables in $\dot{\U}^+_3$}\label{sec4}

In this section we compute the indecomposable objects of $\dot{\U}^+_3$. We show that these are exactly the objects that categorify the elements of the canonical basis set $B$. Furthermore, we 
decompose an arbitrary object as a direct sum of these indecomposables, by categorifying (\ref{glavna}). 

%In such a way we obtain a direct, diagrammatical proof of the result of Varagnolo and Vasserot \cite{vv}. By using geometric %tools, they showed that the indecomposable projective objects categorify the canonical basis, in type $A$.

\begin{theorem}\label{te3}
The objects $\E_1^{(a)}\E_2^{(b)}\E_1^{(c)}\{t\}$ and $\E_2^{(a)}\E_1^{(b)}\E_2^{(c)}\{t\}$, for nonnegative $a,b,c$ and $t\in \Z$, with 
$b\ge a+c$, are indecomposable in $\dot{\U}^+_3$. Neither two of them are isomorphic, except $\E_1^{(a)}\E_2^{(a+c)}\E_1^{(c)}\{t\}\cong
\E_2^{(c)}\E_1^{(a+c)}\E_2^{(a)}\{t\}$, for $a,c\ge 0$, $t\in\Z$.
\end{theorem}

\textbf{Proof:}

First of all, note that it is enough to prove the claim for the objects with no shifts. 
We will show that $\E_1^{(a)}\E_2^{(b)}\E_1^{(c)}$, with $b\ge a+c$, are indecomposable by showing that the graded rank of their endomorphism rings satisfy:
\begin{equation}\label{rang1}
{\rk}_q \Hom\nolimits_{\dot{\U}^*}(\E_1^{(a)}\E_2^{(b)}\E_1^{(c)},\E_1^{(a)}\E_2^{(b)}\E_1^{(c)})\in 1+q\N[q].
\end{equation}
Here $\dot{\U}^*$ is the category with the same objects as $\dot{\U}_3^+$, while morphism between two objects can have arbitary degree. More precisely, for any two objects $x,y\in \dot{\U}^*$, we have $\Hom\nolimits_{\dot{\U}^*}(x,y)=\oplus_{s\in\Z} \Hom\nolimits_{\dot{\U}_3^+}(x\{s\},y)$ (see e.g. \cite{sl2}).

Any diagram from $\E_1^{(a)}\E_2^{(b)}\E_1^{(c)}$ to itself can be represented in a ``standard" form. Namely, we split thick edges into thin ones, and it is well-known (see e.g. \cite{sl2},\cite{kl3}) that any such diagram can be reduced to a diagram where any two thin lines intersect at most once, and all dots from one line are only at one segment. By regrouping the thin edges back into the thick ones, in our case, the only nonzero diagrams are of the form

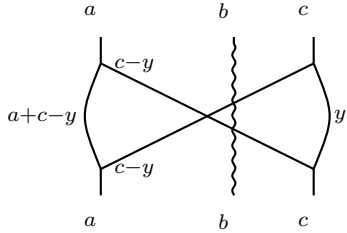
\begin{center}
\begin{tikzpicture}[scale=0.7]
 %\filldraw[black] (4,0) circle (2pt)
%				 (6.5,0) circle (2pt);
				 %(1,-3) circle (2pt)
				 %(1,3) circle (2pt)
				 %(-1,0) circle (2pt)
				 %(1,0) circle (2pt);
 
\draw [thick] (-2,-1.5)--(-2,-1);
\draw [thick] (2,-1.5)--(2,-1);
\draw [thick] (2,1)--(2,1.5);
\draw [thick] (-2,1)--(-2,1.5);

\draw [thick] (-2,1)--(2,-1);
\draw [thick] (-2,-1)--(2,1);
\draw [thick] (-2,-1) .. controls (-2.4,0) .. (-2,1);
\draw [thick] (2,-1) .. controls (2.4,0) .. (2,1);
%\draw  (2,-1.25) .. controls (2.4,0) .. (2,1.25);
\draw [thick,snake=snake,segment amplitude=.2mm,segment length=2mm] 
(0.5,-1.5)--(0.5,1.5);
% \draw (1.5,0) node {=};
\draw (-2.2,-2) node {{{$\scriptstyle{a}$}}};
\draw (-2.2,2) node {{{$\scriptstyle{a}$}}};
\draw (0.3,-2) node {{{$\scriptstyle{b}$}}};
\draw (0.3,2) node {{{$\scriptstyle{b}$}}};
\draw (1.8,-2) node {{{$\scriptstyle{c}$}}};
\draw (1.8,2) node {{{$\scriptstyle{c}$}}};
\draw (-3.1,0) node {{{$\scriptstyle{a+c-y}$}}};
\draw (2.5,0) node {{{$\scriptstyle{y}$}}};
\draw (-1.37,1) node {{{$\scriptstyle{c-y}$}}};
\draw (-1.37,-1) node {{{$\scriptstyle{c-y}$}}};

% \draw (1.2,-3) node {b};
%\draw (2.5,0) node {${{{\sum}_{\alpha\in F}}}$};
%\draw  [->] (4,-3)--(4,3);
%\draw  [->] (6.5,-3)--(6.5,3);

%\draw (3.5,0) node {$\pi_{\alpha}$};
% \draw (7,0) node {$\pi_{\beta}$};
\end{tikzpicture}
\end{center}
possibly with some dots on it. Thus, the lowest degree diagrams are the dotless ones of the form above, where $y\le c$ varies.

Now, for $y=c$ we have the identity - the degree zero map. If $y<c$, the degree of this dotless diagram is 
$$-2(c-y)(a-c+y)-2(c-y)y-2(c-y)^2+2b(c-y)=2(c-y)(b-a-y)>0,$$
for $b\ge a+c$, and so we have proved (\ref{rang1}). Thus, we have that $\E_1^{(a)}\E_2^{(b)}\E_1^{(c)}$ is indecomposable for $b\ge a+c$, as wanted. In the same way, we have that $\E_2^{(a)}\E_1^{(b)}\E_2^{(c)}$ is indecomposable for $b\ge a+c$, as well.\\

In order to prove that two indecomposables of the form $\E_1^{(a)}\E_2^{(b)}\E_1^{(c)}$ and $\E_1^{(p)}\E_2^{(k)}\E_1^{(r)}$, with $(a,b,c)\ne (p,k,r)$ are not isomorphic, we will show that
\begin{equation*}
{\rk}_q \Hom\nolimits_{\dot{\U}^*}(\E_1^{(a)}\E_2^{(b)}\E_1^{(c)},\E_1^{(p)}\E_2^{(k)}\E_1^{(r)})\in q\N[q].
\end{equation*}
This $\Hom$-space can be nonempty only when $k=b$ and $p+r=a+c$, and so we are left with proving:
\begin{equation}\label{rang2}
{\rk}_q \Hom\nolimits_{\dot{\U}^*}(\E_1^{(a+x)}\E_2^{(b)}\E_1^{(c+y)},\E_1^{(a+y)}\E_2^{(b)}\E_1^{(c+x)})\in q\N[q],
\end{equation}
when at least one of $x$ and $y$ is nonzero, with $b\ge a+c+x+y$. Again, a general (dotless) diagram from the last $\Hom$-space has the following form
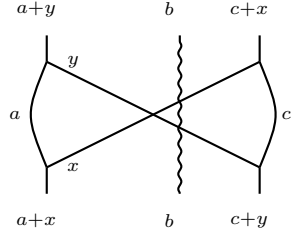
\begin{center}
\begin{tikzpicture}[scale=0.7]
 %\filldraw[black] (4,0) circle (2pt)
%				 (6.5,0) circle (2pt);
				 %(1,-3) circle (2pt)
				 %(1,3) circle (2pt)
				 %(-1,0) circle (2pt)
				 %(1,0) circle (2pt);
 
\draw  [thick] (-2,-1.5)--(-2,-1);
\draw  [thick](2,-1.5)--(2,-1);
\draw  [thick](2,1)--(2,1.5);
\draw  [thick](-2,1)--(-2,1.5);

\draw  [thick](-2,1)--(2,-1);
\draw  [thick](-2,-1)--(2,1);
\draw  [thick](-2,-1) .. controls (-2.4,0) .. (-2,1);
\draw  [thick](2,-1) .. controls (2.4,0) .. (2,1);
%\draw  (2,-1.25) .. controls (2.4,0) .. (2,1.25);
\draw [thick,snake=snake,segment amplitude=.2mm,segment length=2mm] 
(0.5,-1.5)--(0.5,1.5);
% \draw (1.5,0) node {=};
\draw (-2.2,-2) node {\it{\small{$\scriptstyle a+x$}}};
\draw (-2.2,2) node {\it{\small{$\scriptstyle a+y$}}};
\draw (0.3,-2) node {\it{\small{$\scriptstyle b$}}};
\draw (0.3,2) node {\it{\small{$\scriptstyle b$}}};
\draw (1.8,-2) node {\it{\small{$\scriptstyle c+y$}}};
\draw (1.8,2) node {\it{\small{$\scriptstyle c+x$}}};
\draw (-2.6,0) node {\it{\small{$\scriptstyle a$}}};
\draw (2.5,0) node {\it{\small{$\scriptstyle c$}}};
\draw (-1.5,1) node {\it{\small{$\scriptstyle y$}}};
\draw (-1.5,-1) node {\it{\small{$\scriptstyle x$}}};

% \draw (1.2,-3) node {b};
%\draw (2.5,0) node {${{{\sum}_{\alpha\in F}}}$};
%\draw  [->] (4,-3)--(4,3);
%\draw  [->] (6.5,-3)--(6.5,3);

%\draw (3.5,0) node {$\pi_{\alpha}$};
% \draw (7,0) node {$\pi_{\beta}$};
\end{tikzpicture}
\end{center}
The degree of such dotless diagram is equal to
\begin{eqnarray*}&&-ax-ay-cx-cy-2xy+b(x+y)\ge\\
&&\quad\quad\ge -(a+c)(x+y)-2xy+(a+c)(x+y)+(x+y)^2=x^2+y^2>0,
\end{eqnarray*}
as wanted.\\

Finally, to prove that two indecomposables of the form $\E_1^{(a)}\E_2^{(b)}\E_1^{(c)}$ and $\E_2^{(p)}\E_1^{(k)}\E_2^{(r)}$ are not isomorphic, again we compute
\begin{equation*}
{\rk}_q \Hom\nolimits_{\dot{\U}^*}(\E_1^{(a)}\E_2^{(b)}\E_1^{(c)},\E_2^{(p)}\E_1^{(k)}\E_2^{(r)}),
\end{equation*}
for $b\ge a+c$ and $k\ge p+r$. In order this $\Hom$-space to be non-empty, we must have $b=p+r$ and $k=a+c$, and so $b=k=p+r=a+c$. Hence, we are left with computing
\begin{equation}\label{rang4}
{\rk}_q \Hom\nolimits_{\dot{\U}^*}(\E_1^{(a)}\E_2^{(a+c)}\E_1^{(c)},\E_2^{(p)}\E_1^{(a+c)}\E_2^{(r)}),
\end{equation}
with $p+r=a+c$. Again, a general dotless morphism from this $\Hom$-space has the following form 

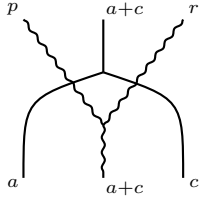
\begin{center}
\begin{tikzpicture}[scale=0.7]
 %\filldraw[black] (4,0) circle (2pt)
%				 (6.5,0) circle (2pt);
				 %(1,-3) circle (2pt)
				 %(1,3) circle (2pt)
				 %(-1,0) circle (2pt)
				 %(1,0) circle (2pt);
 
\draw  [thick] (-1.5,-1.5) .. controls (-1.5,0) .. (0,0.5);
\draw  [thick] (1.5,-1.5) .. controls (1.5,0) .. (0,0.5);
\draw  [thick] (0,0.5)--(0,1.5);

%\draw  (2,-1.5)--(2,-1);
%\draw  (2,1)--(2,1.5);
%\draw  (-2,1)--(-2,1.5);

%\draw  (-2,1)--(2,-1);
%\draw  (-2,-1)--(2,1);
%\draw  (-2,-1) .. controls (-2.4,0) .. (-2,1);
%\draw  (2,-1) .. controls (2.4,0) .. (2,1);
%\draw  (2,-1.25) .. controls (2.4,0) .. (2,1.25);
\draw [thick,snake=snake,segment amplitude=.2mm,segment length=2mm] 
(-1.5,1.5)--(0,-0.5);
\draw [thick,snake=snake,segment amplitude=.2mm,segment length=2mm] 
(1.5,1.5)--(0,-0.5);
\draw [thick,snake=snake,segment amplitude=.2mm,segment length=2mm] 
(0,-0.5)--(0,-1.5);

% \draw (1.5,0) node {=};
\draw (0.4,1.7) node {\it{\small{$\scriptstyle a+c$}}};
\draw (-1.7,-1.6) node {\it{\small{$\scriptstyle a$}}};
\draw (1.7,-1.6) node {\it{\small{$\scriptstyle c$}}};
\draw (0.4,-1.7) node {\it{\small{$\scriptstyle a+c$}}};
\draw (-1.7,1.7) node {\it{\small{$\scriptstyle p$}}};
\draw (1.7,1.7) node {\it{\small{$\scriptstyle r$}}};
%\draw (-2.6,0) node {\it{\small{a}}};
%\draw (2.5,0) node {\it{\small{c}}};
%\draw (-1.5,1) node {\it{\small{y}}};
%\draw (-1.5,-1) node {\it{\small{x}}};

% \draw (1.2,-3) node {b};
%\draw (2.5,0) node {${{{\sum}_{\alpha\in F}}}$};
%\draw  [->] (4,-3)--(4,3);
%\draw  [->] (6.5,-3)--(6.5,3);

%\draw (3.5,0) node {$\pi_{\alpha}$};
% \draw (7,0) node {$\pi_{\beta}$};

 \end{tikzpicture}
 \end{center}

The degree of the diagram from above is equal to
$$-ac-pr+ap+cr=(a-r)(p-c)=(p-c)^2>0,$$
for $p\neq c$. Thus, for $p\ne c$, the graded rank (\ref{rang4}) is from $q\N[q]$, and so we have that the two indecomposable objects $\E_1^{(a)}\E_2^{(b)}\E_1^{(c)}$ and $\E_2^{(p)}\E_1^{(k)}\E_2^{(r)}$ are not isomorphic. 

The only possibility left for the isomorphism is between $\E_1^{(a)}\E_2^{(a+c)}\E_1^{(c)}$ and $\E_2^{(c)}\E_1^{(a+c)}\E_2^{(a)}$, and below (Theorem \ref{gldect}) we show that they are indeed isomorphic. \kraj

\subsection{Decomposition of an arbitrary object of $\dot{\U}^+_3$}

In this section we decompose an arbitrary object from $\dot{\U}^+_3$ as a direct sum of the indecomposable ones from the previous section. We do this by categorifing the Proposition \ref{prgl}, i.e. by proving the following
\begin{theorem}\label{gldect}
For $b\le a+c$, we have the following canonical decomposition
\begin{equation}\label{gldec}
\E_1^{(a)}\E_2^{(b)}\E_1^{(c)}\cong \bigoplus_{\scriptsize{\begin{array}{c}p+r=b\\p\le c\\ r\le a\end{array}}} \bigoplus_{\alpha\in P(c-p,a-r)} \E_2^{(p)}\E_1^{(a+c)}\E_2^{(r)} \{2|\alpha|-(c-p)(a-r)\}.
\end{equation} 
\end{theorem}
Note that this indeed categorifies (\ref{glavna}), since 
$$\left[\begin{array}{c} a+c-b \\ c-p \end{array}\right]=\left[\begin{array}{c} a-r+c-p \\ c-p \end{array}\right]=\sum_{\alpha\in P(c-p,a-r)} q^{2|\alpha|-(c-p)(a-r)}.$$

\begin{remark}
The decomposition (\ref{gldec}) is also valid in general, i.e. if we replace $\E_1$ and $\E_2$ by $\E_r$ and $\E_s$, respectively, with $r\cdot s=-1$.
\end{remark}

Also, in the same way as in the Theorem \ref{mon} (by changing sums with direct sums), the decomposition (\ref{gldec}), together with the decomposition (see [\cite{thick}, Theorem 5.1])  
\[
\E_i^{(a)}\E_i^{(b)}\cong \bigoplus_{\alpha\in P(a,b)} \E_i^{(a+b)}\{2|\alpha|-ab\},
\]
that categorifies (\ref{eiab}), implies
the decomposition of an arbitrary object as a direct sum of the indecomposable ones from the previous section. Thus, the 
decomposition (\ref{gldec}) gives:

\begin{theorem}
The set of indecomposable objects of $\dot{\U}^+_3$ is the following set $\B$:
\begin{equation}
\B=\{\E_1^{(a)}\E_2^{(b)}\E_1^{(c)}\{t\},\,\,\,\E_2^{(a)}\E_1^{(b)}\E_2^{(c)}\{t\}| \quad b\ge a+c,\,\,\, a,b,c\ge 0,\,\,t\in\Z\}. 
\end{equation}
No two elements from $\B$ are isomorphic, except 
\[
\E_1^{(a)}\E_2^{(a+c)}\E_1^{(c)}\{t\}\cong
\E_2^{(c)}\E_1^{(a+c)}\E_2^{(a)}\{t\},\quad a,c\ge 0,\,\,t\in\Z.
\]
An arbitrary object of $\dot{\U}^+_3$ can be decomposed as a direct sum of the elements from $\B$.
\end{theorem}

In this way we have obtained a bijection between the canonical basis $B$ of $U_q^+(\sl_3)$ and the indecomposable objects from $\B$
with no shifts, given by:
\begin{eqnarray*}
E_1^{(a)}E_2^{(b)}E_1^{(c)}&\mapsto & \E_1^{(a)}\E_2^{(b)}\E_1^{(c)}, \quad b\ge a+c,\\
E_2^{(a)}E_1^{(b)}E_2^{(c)}&\mapsto & \E_2^{(a)}\E_1^{(b)}\E_2^{(c)}, \quad b\ge a+c,\\
E_1^{(a)}E_2^{(a+c)}E_1^{(c)}=E_2^{(c)}E_1^{(a+c)}E_2^{(a)}&\mapsto &\E_1^{(a)}\E_2^{(a+c)}\E_1^{(c)}\cong\E_2^{(c)}\E_1^{(a+c)}\E_2^{(a)}.
\end{eqnarray*}
\vskip 0.5cm

So, we are left with proving the decomposition (\ref{gldec}) (i.e. Theorem \ref{gldect}) and that is done in the rest of this section. \\

\subsection{Proof of Theorem \ref{gldect}}

First of all, the conditions $p+r=b$, $p\le c$ and $r\le a$ together, are equivalent to $\max(0,b-a)\le p \le \min(b,c)$, with $r=b-p$. Then, for every non-negative integer $p$ with $\max(0,b-a)\le p \le \min(b,c)$, and partition $\alpha\in P(c-p,a-r)$, where $r=b-p$, we define the following 2-morphisms:

\begin{equation}\label{lamb}
\lambda_{\alpha}^p: \E_2^{(p)}\E_1^{(a+c)}\E_2^{(r)}\{2|\alpha|-(c-p)(a-r)\}\longrightarrow \E_1^{(a)}\E_2^{(b)}\E_1^{(c)},
\end{equation}

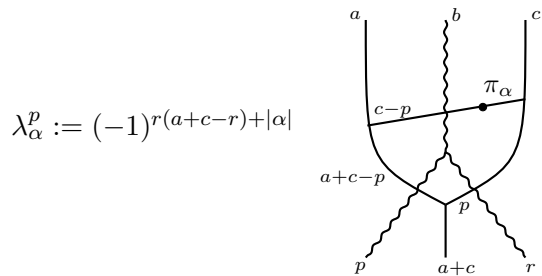
\begin{center}
\begin{tikzpicture}[scale=0.7]
 \filldraw[black] (0.7,1.35) circle (2pt);
%				 (6.5,0) circle (2pt);
				 %(1,-3) circle (2pt)
				 %(1,3) circle (2pt)
				 %(-1,0) circle (2pt)
				 %(1,0) circle (2pt);
\draw (-5.5,1) node {$\lambda^p_{\alpha}:=(-1)^{r(a+c-r)+|\alpha|}$}; 
\draw  [thick](-1.5,3) .. controls (-1.5,0.3) .. (0,-0.5);
\draw  [thick](1.5,3) .. controls (1.5,0.3) .. (0,-0.5);
\draw  [thick](0,-0.5)--(0,-1.5);
\draw  [thick](-1.47,1)--(1.5,1.5);

\draw [thick,snake=snake,segment amplitude=.2mm,segment length=2mm] 
(0,3)--(0,0.5);
\draw [thick,snake=snake,segment amplitude=.2mm,segment length=2mm] 
(0,0.5)--(-1.5,-1.5);
\draw [thick,snake=snake,segment amplitude=.2mm,segment length=2mm] 
(0,0.5)--(1.5,-1.5);

% \draw (1.5,0) node {=};
\draw (-1.6,-1.7) node {\it{\small{$\scriptstyle p$}}};
\draw (0.2,-1.7) node {\it{\small{$\scriptstyle a+c$}}};
\draw (1.6,-1.7) node {\it{\small{$\scriptstyle r$}}};
\draw (0.35,-0.55) node {\it{\small{$\scriptstyle p$}}};
\draw (-1.75,0) node {\it{\small{$\scriptstyle a+c-p$}}};
\draw (-1,1.3) node {\it{\small{$\scriptstyle c-p$}}};
\draw (1,1.7) node {\it{\small{$\pi_{\alpha}$}}};
\draw (-1.7,3.1) node {\it{\small{$\scriptstyle a$}}};
\draw (0.2,3.1) node {\it{\small{$\scriptstyle b$}}};
\draw (1.7,3.1) node {\it{\small{$\scriptstyle c$}}};

%\draw (-1.7,-1.6) node {\it{\small{a}}};
%\draw (1.7,-1.6) node {\it{\small{c}}};
%\draw (0.4,-1.7) node {\it{\small{a+c}}};
%\draw (-1.7,1.7) node {\it{\small{p}}};
%\draw (1.7,1.7) node {\it{\small{r}}};

 \end{tikzpicture}
\end{center}

\begin{equation}\label{sig}
\sigma_{\alpha}^p: \E_1^{(a)}\E_2^{(b)}\E_1^{(c)}\longrightarrow \E_2^{(p)}\E_1^{(a+c)}\E_2^{(r)}\{2|\alpha|-(c-p)(a-r)\},
\end{equation} 

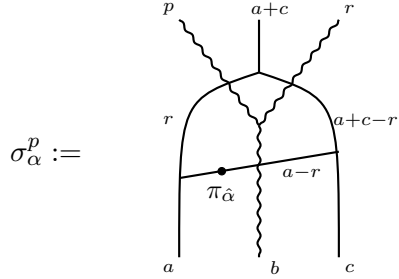
\begin{center}
\begin{tikzpicture}[scale=0.7]
 \filldraw[black] (-0.7,-1.37) circle (2pt);
%				 (6.5,0) circle (2pt);
				 %(1,-3) circle (2pt)
				 %(1,3) circle (2pt)
				 %(-1,0) circle (2pt)
				 %(1,0) circle (2pt);

\draw (-4,-1) node {$\sigma^p_{\alpha}:=$}; 
\draw  [thick](0,0.5) .. controls (1.5,0) .. (1.5,-3);
\draw  [thick](0,0.5) .. controls (-1.5,0) .. (-1.5,-3);
\draw  [thick](0,1.5)--(0,0.5);
\draw  [thick](-1.5,-1.5)--(1.5,-1);

\draw [thick,snake=snake,segment amplitude=.2mm,segment length=2mm] 
(-1.5,1.5)--(0,-0.5);
\draw [thick,snake=snake,segment amplitude=.2mm,segment length=2mm] 
(1.5,1.5)--(0,-0.5);
\draw [thick,snake=snake,segment amplitude=.2mm,segment length=2mm] 
(0,-0.5)--(0,-3);

% \draw (1.5,0) node {=};
\draw (-1.7,-3.2) node {\it{\small{$\scriptstyle a$}}};
\draw (0.3,-3.2) node {\it{\small{$\scriptstyle b$}}};
\draw (1.7,-3.2) node {\it{\small{$\scriptstyle c$}}};
\draw (-1.7,1.7) node {\it{\small{$\scriptstyle p$}}};
\draw (0.2,1.7) node {\it{\small{$\scriptstyle a+c$}}};
\draw (1.7,1.7) node {\it{\small{$\scriptstyle r$}}};
\draw (-1.7,-0.4) node {\it{\small{$\scriptstyle r$}}};
\draw (0.8,-1.4) node {\it{\small{$\scriptstyle a-r$}}};
\draw (2,-0.4) node {\it{\small{$\scriptstyle a+c-r$}}};
\draw (-0.7,-1.8) node {\it{\small{$\pi_{\hat{\alpha}}$}}};
 \end{tikzpicture}
\end{center}

\begin{equation}\label{eovi}
e^p_{\alpha}:=\lambda_{\alpha}^p \sigma_{\alpha}^p: \E_1^{(a)}\E_2^{(b)}\E_1^{(c)} \longrightarrow \E_1^{(a)}\E_2^{(b)}\E_1^{(c)}.
\end{equation}\\

The following lemma is the key result:

\begin{lemma}\label{gl}
Let $b\le a+c$. Let $\max(0,b-a)\le p,p' \le \min(b,c)$, $\alpha\in P(c-p,a-r)$ and $\alpha'\in P(c-p',a-r')$, where $r=b-p$ 
and $r'=b-p'$. Then:
\begin{equation}
\sigma_{\alpha'}^{p'} \lambda_{\alpha}^p = \delta_{p,p'} \delta_{\alpha,\alpha'} {\Id}_{\E_2^{(p)}\E_1^{(a+c)}\E_2^{(r)}}.
\end{equation}
\end{lemma}

The last lemma implies the theorem below, which in turn proves the wanted decomposition (\ref{gldec}).

\begin{theorem}
The collection $\{e_{\alpha}^p\}$ is a collection of mutually orthogonal idempotents. 
\end{theorem}
%Thus, all $\{e_{\alpha}^p\}$ are idempotents on $\E_1^{(a)}\E_2^{(b)}\E_1^{(c)}$ that factor through $\E_2^{(p)}\E_1^{(a+c)}\E_2^{(r)}\{2|\alpha|-(c-p)(a-r)\}$.\\

\textbf{Proof of Lemma \ref{gl}:}

%\begin{equation}
%\pi_{\alpha}(\underline{x})=\sum_{\beta,\gamma\in P(c-p)} (-1)^{|\gamma|} c_{\beta,\gamma}^{\alpha} \pi_{\beta}(\underline{x},
%\underline{y}) \pi_{\bar{\gamma}}(\underline{y}).
%\end{equation}

In pictures, the statement of the lemma above is the following:

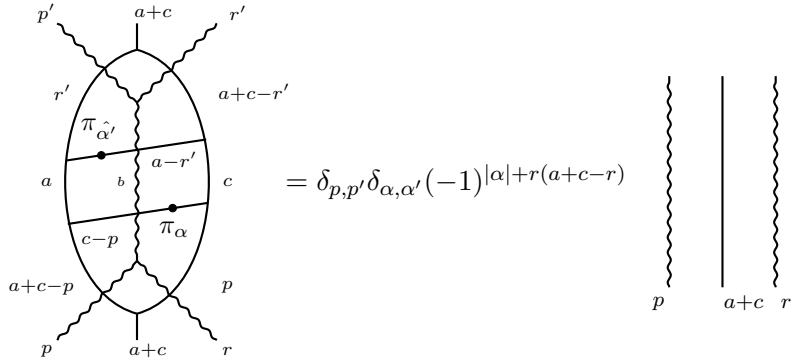
\begin{center}
\begin{tikzpicture}[scale=0.7]
 \filldraw[black] (-0.67,0.5) circle (2pt)
				 (0.67,-0.5) circle (2pt);
				 %(1,-3) circle (2pt)
				 %(1,3) circle (2pt)
				 %(-1,0) circle (2pt)
				 %(1,0) circle (2pt);
 
\draw  [thick](0,2.5) .. controls (-1.8,2) and (-1.8,-2) .. (0,-2.5);
\draw  [thick](0,2.5) .. controls (1.8,2) and (1.8,-2) .. (0,-2.5);
\draw  [thick](0,3)--(0,2.5);
\draw  [thick](0,-3)--(0,-2.5);
\draw  [thick](-1.34,0.4)--(1.3,0.8);
\draw  [thick](-1.3,-0.8)--(1.34,-0.4);

\draw [thick,snake=snake,segment amplitude=.2mm,segment length=2mm] 
(-1.5,3)--(0,1.5);
\draw [thick,snake=snake,segment amplitude=.2mm,segment length=2mm] 
(1.5,3)--(0,1.5);
\draw [thick,snake=snake,segment amplitude=.2mm,segment length=2mm] 
(0,1.5)--(0,-1.5);
\draw [thick,snake=snake,segment amplitude=.2mm,segment length=2mm] 
(-1.5,-3)--(0,-1.5);
\draw [thick,snake=snake,segment amplitude=.2mm,segment length=2mm] 
(1.5,-3)--(0,-1.5);

% \draw (1.5,0) node {=};
\draw (-1.7,3.2) node {\it{\small{$\scriptstyle p'$}}};
\draw (0.3,3.2) node {\it{\small{$\scriptstyle a+c$}}};
\draw (1.9,3.2) node {\it{\small{$\scriptstyle r'$}}};
\draw (-1.7,-3.2) node {\it{\small{$\scriptstyle p$}}};
\draw (0.2,-3.2) node {\it{\small{$\scriptstyle a+c$}}};
\draw (1.7,-3.2) node {\it{\small{$\scriptstyle r$}}};

\draw (-0.7,1) node {\it{\small{$\pi_{\hat{\alpha'}}$}}};
\draw (0.7,0.4) node {\it{\small{$\scriptstyle a-r'$}}};

\draw (-1.4,1.7) node {\it{\small{$\scriptstyle r'$}}};
\draw (2.2,1.7) node {\it{\small{$\scriptstyle a+c-r'$}}};

\draw (-0.3,0) node {\it{\tiny{$\scriptstyle{b}$}}};

\draw (-1.7,0) node {\it{\small{$\scriptstyle a$}}};

\draw (1.7,0) node {\it{\small{$\scriptstyle c$}}};

\draw (0.7,-0.9) node {\it{\small{$\pi_{{\alpha}}$}}};
\draw (-0.7,-1.1) node {\it{\small{$\scriptstyle c-p$}}};

\draw (-1.8,-2) node {\it{\small{$\scriptstyle a+c-p$}}};
\draw (1.7,-2) node {\it{\small{$\scriptstyle p$}}};

\draw (6,0) node {$= \delta_{p,p'}\delta_{\alpha,\alpha'}(-1)^{|\alpha|+r(a+c-r)}$};
\draw [thick,snake=snake,segment amplitude=.2mm,segment length=2mm] 
(10,-2)--(10,2);
\draw [thick](11,-2)--(11,2);
\draw [thick,snake=snake,segment amplitude=.2mm,segment length=2mm] (12,-2)--(12,2);
\draw (9.8,-2.3) node {$\scriptstyle{p}$};
\draw (11.4,-2.3) node {$\scriptstyle{a+c}$};
\draw (12.2,-2.3) node {$\scriptstyle{r}$};

\end{tikzpicture}
\end{center}

We shall prove the formula by simplifying the diagram on the left hand side: 

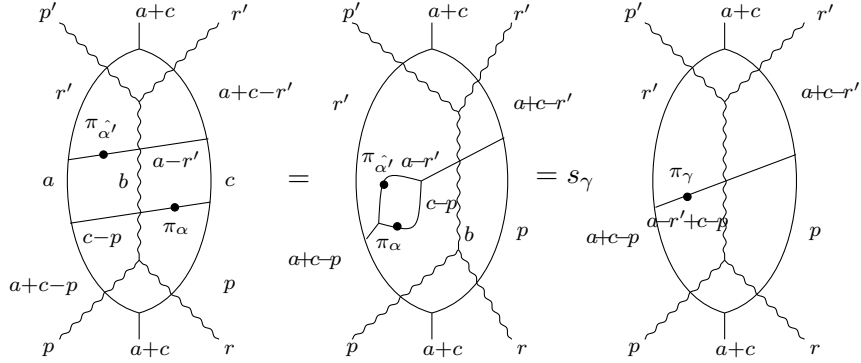
\begin{center}
\begin{tikzpicture}[scale=0.7]

 \filldraw[black]  (-0.67,0.5) circle (2pt)
				 (0.67,-0.5) circle (2pt);
				 %(1,-3) circle (2pt)
				 %(1,3) circle (2pt)
				 %(-1,0) circle (2pt)
				 %(1,0) circle (2pt);
 
\draw  (0,2.5) .. controls (-1.8,2) and (-1.8,-2) .. (0,-2.5);
\draw  (0,2.5) .. controls (1.8,2) and (1.8,-2) .. (0,-2.5);
\draw  (0,3)--(0,2.5);
\draw  (0,-3)--(0,-2.5);
\draw  (-1.34,0.4)--(1.3,0.8);
\draw  (-1.3,-0.8)--(1.34,-0.4);

\draw [snake=snake,segment amplitude=.2mm,segment length=2mm] 
(-1.5,3)--(0,1.5);
\draw [snake=snake,segment amplitude=.2mm,segment length=2mm] 
(1.5,3)--(0,1.5);
\draw [snake=snake,segment amplitude=.2mm,segment length=2mm] 
(0,1.5)--(0,-1.5);
\draw [snake=snake,segment amplitude=.2mm,segment length=2mm] 
(-1.5,-3)--(0,-1.5);
\draw [snake=snake,segment amplitude=.2mm,segment length=2mm] 
(1.5,-3)--(0,-1.5);

% \draw (1.5,0) node {=};
\draw (-1.7,3.2) node {$\scriptstyle{p'}$};
\draw (0.3,3.2) node {$\scriptstyle{a+c}$};
\draw (1.9,3.2) node {$\scriptstyle{r'}$};
\draw (-1.7,-3.2) node {$\scriptstyle{p}$};
\draw (0.2,-3.2) node {$\scriptstyle{a+c}$};
\draw (1.7,-3.2) node {$\scriptstyle{r}$};

\draw (-0.7,1) node {$\scriptstyle{ \pi_{\hat{\alpha'}}}$};
\draw (0.7,0.4) node {$\scriptstyle{ a-r'}$};

\draw (-1.4,1.7) node {$\scriptstyle{r'}$};
\draw (2.2,1.7) node {$\scriptstyle{a+c-r'}$};

\draw (-0.3,0) node {$\scriptstyle{b}$};

\draw (-1.7,0) node {$\scriptstyle{a}$};

\draw (1.7,0) node {$\scriptstyle{c}$};

\draw (0.7,-0.9) node {$\scriptstyle{\pi_{{\alpha}}}$};
\draw (-0.7,-1.1) node {$\scriptstyle{c-p}$};

\draw (-1.8,-2) node {$\scriptstyle{a+c-p}$};
\draw (1.7,-2) node {$\scriptstyle{p}$};

\draw (3,0) node {$\displaystyle =$};

 \filldraw[black] (4.6,-0.07) circle (2pt)
				 (4.85,-0.86) circle (2pt);
				 %(1,-3) circle (2pt)
				 %(1,3) circle (2pt)
				 %(-1,0) circle (2pt)
				 %(1,0) circle (2pt);
 
\draw  (5.5,2.5) .. controls (3.6,2) and (3.6,-2) .. (5.5,-2.5);
\draw  (5.5,2.5) .. controls (7.4,2) and (7.4,-2) .. (5.5,-2.5);
\draw  (5.5,3)--(5.5,2.5);
\draw  (5.5,-3)--(5.5,-2.5);

\draw  (4.27,-1.1)--(4.5,-0.8);
\draw  (4.5,-0.8)..controls (4.55,0.2)..(5.3,0);
\draw  (4.5,-0.8)..controls (5.25,-1)..(5.3,0);
\draw (5.3,0)--(6.86,0.82);

\draw [snake=snake,segment amplitude=.2mm,segment length=2mm] 
(4,3)--(6,1.3);
\draw [snake=snake,segment amplitude=.2mm,segment length=2mm] 
(7,3)--(6,1.3);
\draw [snake=snake,segment amplitude=.2mm,segment length=2mm] 
(6,1.3)--(6,-1.3);
\draw [snake=snake,segment amplitude=.2mm,segment length=2mm] 
(4,-3)--(6,-1.3);
\draw [snake=snake,segment amplitude=.2mm,segment length=2mm] 
(7,-3)--(6,-1.3);

% \draw (1.5,0) node {=};
\draw (4,3.2) node {{$\scriptstyle{p'}$}};
\draw (5.8,3.2) node {$\scriptstyle{a+c}$};
\draw (7.2,3.2) node {$\scriptstyle{r'}$};

\draw (4,-3.2) node {$\scriptstyle{p}$};
\draw (5.7,-3.2) node {$\scriptstyle{a+c}$};
\draw (7.2,-3.2) node {$\scriptstyle{r}$};

\draw (3.8,1.5) node {$\scriptstyle{r'}$};
\draw (3.3,-1.5) node {$\scriptstyle{a\!+\!c\!-\!p}$};

\draw (7.6,1.5) node {$\scriptstyle{a\!+\!c\!-\!r'}$};
\draw (7.2,-1) node {$\scriptstyle{p}$};

\draw (6.2,-1) node {$\scriptstyle{b}$};

\draw (5.3,0.4) node {$\scriptstyle{a\!-\!r'}$};
\draw (5.68,-0.45) node {$\scriptstyle{c\!-\!p}$};

\draw (4.5,0.35) node {$\scriptstyle{\pi_{\hat{\alpha'}}}$};
%\draw (0.7,0.4) node {\it{\small{$a-r'$}}};

%\draw (-0.4,0) node {\it{\tiny{p+r}}};

%\draw (-1.7,0) node {\it{\small{a}}};

%\draw (1.7,0) node {\it{\small{c}}};

\draw (4.7,-1.2) node {$\scriptstyle{\pi_{{\alpha}}}$};
%\draw (-0.7,-1.1) node {\it{\small{c-p}}};

%\draw (-1.8,-2) node {\it{\small{a+c-p}}};
%\draw (1.7,-2) node {\it{\small{p}}};
\draw (8,0) node {$\displaystyle = s_{\gamma}$};

 \filldraw[black]  (10.3,-0.3) circle (2pt);
%				 (0.67,-0.5) circle (2pt);
				 %(1,-3) circle (2pt)
				 %(1,3) circle (2pt)
				 %(-1,0) circle (2pt)
				 %(1,0) circle (2pt);
 
\draw  (11,2.5) .. controls (9.2,2) and (9.2,-2) .. (11,-2.5);
\draw  (11,2.5) .. controls (12.8,2) and (12.8,-2) .. (11,-2.5);
\draw  (11,3)--(11,2.5);
\draw  (11,-3)--(11,-2.5);
%\draw  (-1.34,0.4)--(1.3,0.8);
\draw  (9.7,-0.5)--(12.34,0.5);

\draw [snake=snake,segment amplitude=.2mm,segment length=2mm] 
(9.5,3)--(11,1.5);
\draw [snake=snake,segment amplitude=.2mm,segment length=2mm] 
(12.5,3)--(11,1.5);
\draw [snake=snake,segment amplitude=.2mm,segment length=2mm] 
(11,1.5)--(11,-1.5);
\draw [snake=snake,segment amplitude=.2mm,segment length=2mm] 
(9.5,-3)--(11,-1.5);
\draw [snake=snake,segment amplitude=.2mm,segment length=2mm] 
(12.5,-3)--(11,-1.5);

% \draw (1.5,0) node {=};
\draw (9.3,3.2) node {$\scriptstyle{p'}$};
\draw (11.3,3.2) node {$\scriptstyle{a+c}$};
\draw (12.9,3.2) node {$\scriptstyle{r'}$};

\draw (9.3,-3.2) node {$\scriptstyle{p}$};
\draw (11.2,-3.2) node {$\scriptstyle{a+c}$};
\draw (12.7,-3.2) node {$\scriptstyle{r}$};

\draw (9.6,1.7) node {$\scriptstyle{r'}$};
\draw (8.9,-1.1) node {$\scriptstyle{a\!+\!c\!-\!p}$};

\draw (13,1.7) node {$\scriptstyle{a\!+\!c\!-\!r'}$};
\draw (12.7,-1) node {$\scriptstyle{p}$};

\draw (10.3,-0.75) node {$\scriptstyle{a\!-\!r'\!+\!c\!-\!p}$};

\draw (10.2,0.1) node {$\scriptstyle{\pi_{\gamma}}$};

%\draw (15,0) node {$\displaystyle =$};
\end{tikzpicture}
\end{center}

\noindent where, by Lemma \ref{vazlem}, $\gamma\in P(a\!-\!r'\!\!+\!c\!-\!p,r'\!\!-\!r)$ and $s_{\gamma}\in\{-1,0,1\}$. Thus, in order the last diagram to be nonzero, we must have $r'\ge r$. Moreover, if $r'=r$, by the second part of Lemma \ref{vazlem}, we must also have $\alpha=\alpha'$ and $s_{\gamma}=(-1)^{|\alpha|}$. 

Now, the last diagram (without the sign $s_{\gamma}$), by applying the Pitchfork lemma, Thick R2 move and Lemma \ref{pomocna}, becomes:

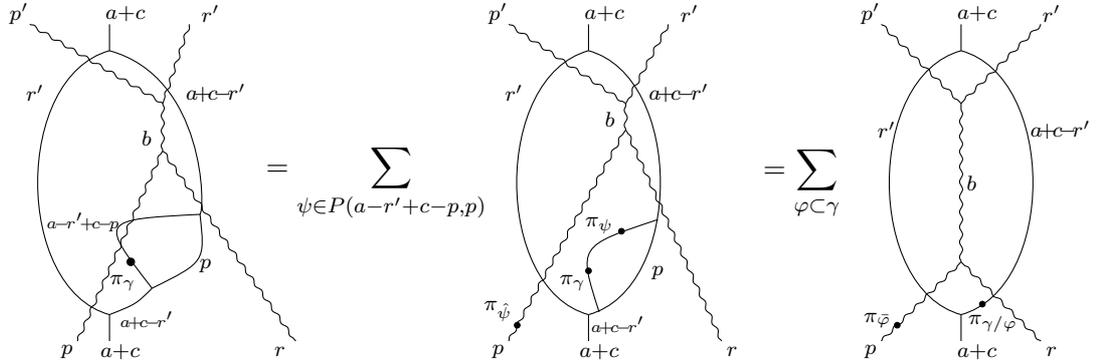
\begin{center}
\begin{tikzpicture}[scale=0.7]

 \filldraw[black]  (11.4,-1.5) circle (2pt);
%				 (0.67,-0.5) circle (2pt);
				 %(1,-3) circle (2pt)
				 %(1,3) circle (2pt)
				 %(-1,0) circle (2pt)
				 %(1,0) circle (2pt);
 
\draw  (11,2.5) .. controls (9.2,2) and (9.2,-2) .. (11,-2.5);
\draw  (11,2.5) .. controls (12.8,2) and (12.8,-0.5) .. (12.7,-0.6);

\draw  (12.7,-0.6) .. controls (12.8,-1.5) .. (11.8,-2);

\draw  (12.7,-0.6) .. controls (10.8,-0.7) .. (11.8,-2);

\draw (11.8,-2)..controls (11.5,-2.3)..(11,-2.5);

\draw  (11,3)--(11,2.5);
\draw  (11,-3)--(11,-2.5);
%\draw  (-1.34,0.4)--(1.3,0.8);
%\draw  (9.7,-0.5)--(12.34,0.5);

\draw [snake=snake,segment amplitude=.2mm,segment length=2mm] 
(9.5,3)--(12,1.5);
\draw [snake=snake,segment amplitude=.2mm,segment length=2mm] 
(12.5,3)--(12,1.5);
\draw [snake=snake,segment amplitude=.2mm,segment length=2mm] 
(12,1.5)--(12,0.6);
\draw [snake=snake,segment amplitude=.2mm,segment length=2mm] 
(10.4,-3)--(12,0.6);
\draw [snake=snake,segment amplitude=.2mm,segment length=2mm] 
(12,0.6)--(14.5,-3);

% \draw (1.5,0) node {=};
\draw (9.3,3.2) node {${\scriptstyle{p'}}$};
\draw (11.3,3.2) node {${\scriptstyle{a+c}}$};
\draw (12.9,3.2) node {${\scriptstyle{r'}}$};

\draw (10.2,-3.2) node {${\scriptstyle{p}}$};
\draw (11.2,-3.2) node {${\scriptstyle{a+c}}$};
\draw (14.2,-3.2) node {${\scriptstyle{r}}$};

\draw (9.6,1.7) node {${\scriptstyle{r'}}$};
%\draw (8.9,-1.1) node {\it{\small{$a\!+\!c\!-\!p$}}};

\draw (11.7,0.85) node {$\scriptstyle b$};

\draw (13,1.7) node {${\scriptstyle{a\!+\!c\!-\!r'}}$};
\draw (12.8,-1.57) node {${\scriptstyle{p}}$};
\draw (11.7,-2.6) node {${\scriptscriptstyle{a\!+\!c\!-\!r'}}$};

\draw (10.5,-0.75) node {{{$\scriptscriptstyle{a\!-\!r'\!+\!c\!-\!p}$}}};

\draw (11.25,-1.9) node {{{$\scriptstyle \pi_{\gamma}$}}};

%\draw (1.7,0) node {\it{\small{c}}};

\draw (16,0) node {$\displaystyle{= \sum_{\psi\in P(a-r'+c-p,p)}}$};

 \filldraw[black]  (18.66,-2.7) circle (1.5pt)
				 (20,-1.67) circle (1.5pt)
				 (20.62,-0.92) circle (1.5pt);
				 %(1,3) circle (2pt)
				 %(-1,0) circle (2pt)
				 %(1,0) circle (2pt);
 
 \draw (18.3,-2.43) node {$\scriptstyle{\pi_{\hat{\psi}}}$};
 \draw (19.7,-1.9) node {$\scriptstyle{\pi_{{\gamma}}}$};
 \draw (20.2,-0.8) node {$\scriptstyle{\pi_{{\psi}}}$};

\draw  (20,2.5) .. controls (18.2,2) and (18.2,-2) .. (20,-2.5);
\draw  (20,2.5) .. controls (21.8,2) and (21.8,-2) .. (20,-2.5);
\draw  (20,3)--(20,2.5);
\draw  (20,-3)--(20,-2.5);
%\draw  (-1.34,0.4)--(1.3,0.8);
%\draw  (18.7,-0.5)--(21.34,0.5);
\draw (20.2,-2.45)..controls (19.8,-1.2)..(21.3,-0.7);

\draw [snake=snake,segment amplitude=.2mm,segment length=2mm] 
(18.5,3)--(20.7,1.5);
\draw [snake=snake,segment amplitude=.2mm,segment length=2mm] 
(21.5,3)--(20.7,1.5);
\draw [snake=snake,segment amplitude=.2mm,segment length=2mm] 
(20.7,1)--(20.7,1.5);
\draw [snake=snake,segment amplitude=.2mm,segment length=2mm] 
(18.5,-3)--(20.7,1);
\draw [snake=snake,segment amplitude=.2mm,segment length=2mm] 
(22.5,-3)--(20.7,1);

% \draw (1.5,0) node {=};
\draw (18.3,3.2) node {{{$\scriptstyle p'$}}};
\draw (20.3,3.2) node {${\scriptstyle {a+c}}$};
\draw (21.9,3.2) node {${{\scriptstyle r'}}$};

\draw (18.3,-3.2) node {${\scriptstyle{p}}$};
\draw (20.2,-3.2) node {${\scriptstyle{a+c}}$};
\draw (22.7,-3.2) node {${\scriptstyle{r}}$};

\draw (18.6,1.7) node {${\scriptstyle{r'}}$};
%\draw (17.9,-1.1) node {\it{\small{$a\!+\!c\!-\!p$}}};

\draw (21.7,1.7) node {{{$ \scriptstyle a\!+\!c\!-\!r'$}}};
\draw (21.3,-1.7) node {{{$\scriptstyle p$}}};

%\draw (19.3,-0.75) node {\it{\tiny{$\scriptstyle{a\!-\!r'\!+\!c\!-\!p}$}}};

%\draw (19.2,0.1) node {{{$\scriptstyle \pi_{\gamma}$}}};

\draw (20.4,1.2) node {{{$\scriptstyle b$}}};
\draw (20.55,-2.7) node {{{$\scriptscriptstyle a\!+\!c\!-\!r'$}}};

\draw (24,0) node {$\displaystyle{= \sum_{\varphi\subset \gamma}}$};

 \filldraw[black]  (25.8,-2.7) circle (1.5pt)
				 (27.4,-2.3) circle (1.5pt);
				 %(20.62,-0.92) circle (1.5pt);
				 %(1,3) circle (2pt)
				 %(-1,0) circle (2pt)
				 %(1,0) circle (2pt);
 
 \draw (25.43,-2.6) node {$\scriptstyle{\pi_{\bar{\varphi}}}$};
 \draw (27.6,-2.67) node {$\scriptstyle{\pi_{{\gamma}/ \varphi}}$};
 %\draw (20.2,-0.8) node {$\scriptstyle{\pi_{{\psi}}}$};

\draw  (27,2.5) .. controls (25.2,2) and (25.2,-2) .. (27,-2.5);
\draw  (27,2.5) .. controls (28.8,2) and (28.8,-2) .. (27,-2.5);
\draw  (27,3)--(27,2.5);
\draw  (27,-3)--(27,-2.5);
%\draw  (-1.34,0.4)--(1.3,0.8);
%\draw  (18.7,-0.5)--(21.34,0.5);
%\draw (28.2,-2.45)..controls (27.8,-1.2)..(29.3,-0.7); 

\draw [snake=snake,segment amplitude=.2mm,segment length=2mm] 
(25.5,3)--(27,1.5);
\draw [snake=snake,segment amplitude=.2mm,segment length=2mm] 
(28.5,3)--(27,1.5);
\draw [snake=snake,segment amplitude=.2mm,segment length=2mm] 
(25.5,-3)--(27,-1.5);
\draw [snake=snake,segment amplitude=.2mm,segment length=2mm] 
(28.5,-3)--(27,-1.5);
\draw [snake=snake,segment amplitude=.2mm,segment length=2mm] 
(27,-1.5)--(27,1.5);

% \draw (1.5,0) node {=};
\draw (25.3,3.2) node {{{$\scriptstyle p'$}}};
\draw (27.3,3.2) node {${\scriptstyle {a+c}}$};
\draw (28.7,3.2) node {${{\scriptstyle r'}}$};

\draw (25.3,-3.2) node {${\scriptstyle{p}}$};
\draw (27.3,-3.2) node {${\scriptstyle{a+c}}$};
\draw (28.7,-3.2) node {${\scriptstyle{r}}$};

\draw (25.6,1) node {${\scriptstyle{r'}}$};
%\draw (17.9,-1.1) node {\it{\small{$a\!+\!c\!-\!p$}}};

\draw (28.87,1) node {{{$ \scriptstyle a\!+\!c\!-\!r'$}}};
\draw (27.2,0) node {{{$\scriptstyle b$}}};

%\draw (19.3,-0.75) node {\it{\tiny{$\scriptstyle{a\!-\!r'\!+\!c\!-\!p}$}}};

%\draw (19.2,0.1) node {{{$\scriptstyle \pi_{\gamma}$}}};

%\draw (20.4,1.2) node {{{$\scriptstyle b$}}};
%\draw (20.55,-2.7) node {{{$\scriptscriptstyle a\!+\!c\!-\!r'$}}};
\end{tikzpicture}
\end{center} 

On the last diagram, we shall apply Opening of a Thick Edge, for the curly line of thickness $b$. We have two possibilities: either $p\ge r'$ or $p<r'$. We shall assume that the first one is satisfied - the other case is done completely analogously.

So, let $p=r'+x$, for some $x\ge 0$. Note that then also $p'=r+x$. Then, the last diagram becomes:

\begin{center}
\begin{tikzpicture}[scale=1.1]
\draw (-3,0) node {$\displaystyle{\sum_{\varphi\subset \gamma}}$};
\filldraw[black]  (-1.4,-1.13) circle (1.5pt)
				 (-0.6,-0.9) circle (1.5pt);
				 %(20.62,-0.92) circle (1.5pt);
				 %(1,3) circle (2pt)
				 %(-1,0) circle (2pt)
				 %(1,0) circle (2pt);
 
 \draw (-1.7,-1.1) node {$\scriptstyle{\pi_{\bar{\varphi}}}$};
 \draw (-0.17,-1.04) node {$\scriptstyle{\pi_{{\gamma}/ \varphi}}$};
 %\draw (20.2,-0.8) node {$\scriptstyle{\pi_{{\psi}}}$};

\draw  (-1,1.25) .. controls (-1.8,1) and (-1.8,-1) .. (-1,-1.25);
\draw  (-1,1.25) .. controls (-0.2,1) and (-0.2,-1) .. (-1,-1.25);
\draw  (-1,1.5)--(-1,1.25);
\draw  (-1,-1.5)--(-1,-1.25);
%\draw  (-1.34,0.4)--(1.3,0.8);
%\draw  (18.7,-0.5)--(21.34,0.5);
%\draw (28.2,-2.45)..controls (27.8,-1.2)..(29.3,-0.7); 

\draw [snake=snake,segment amplitude=.2mm,segment length=2mm] 
(-1.5,1.3)--(0.3,-1.3);
\draw [snake=snake,segment amplitude=.2mm,segment length=2mm] 
(-1.5,-1.3)--(0.3,1.3);
\draw [snake=snake,segment amplitude=.2mm,segment length=2mm] 
(-1.2,-0.85)--(-1.2,0.86);
%\draw [snake=snake,segment amplitude=.2mm,segment length=2mm] 
%(28.5,-3)--(27,-1.5);
%\draw [snake=snake,segment amplitude=.2mm,segment length=2mm] 
%(27,-1.5)--(27,1.5);

% \draw (1.5,0) node {=};
\draw (-1.6,1.5) node {{{$\scriptstyle p'$}}};
\draw (-0.8,1.7) node {${\scriptstyle {a+c}}$};
\draw (0.5,1.5) node {${{\scriptstyle r'}}$};

\draw (-1.6,-1.5) node {${\scriptstyle{p}}$};
\draw (-0.8,-1.7) node {${\scriptstyle{a+c}}$};
\draw (0.5,-1.5) node {${\scriptstyle{r}}$};

\draw (-1.77,0) node {${\scriptstyle{r'}}$};
%\draw (17.9,-1.1) node {\it{\small{$a\!+\!c\!-\!p$}}};

\draw (0.12,0) node {{{$ \scriptstyle a\!+\!c\!-\!r'$}}};
\draw (-1.4,0) node {{{$\scriptstyle x$}}};

\draw (-0.9,0.7) node {{{$\scriptstyle r$}}};

\draw (-0.9,-0.7) node {{{$\scriptstyle r'$}}};

 \filldraw[black]  (5.18,-0.57) circle (1.7pt)
(3.9,-1.1) circle (1.7pt);

\draw (1.5,0.4) node {{\tiny{Pitchfork}}};
\draw (2,0) node {$\displaystyle{=\quad\sum_{\varphi\subset \gamma}}$};

\draw (5,-1.5)--(5,-0.7);
\draw (5,0.7)--(5,1.5);
\draw (5,-0.7)..controls (3.5,0).. (5,0.7);
\draw (5,-0.7)..controls (6,0).. (5,0.7);
\draw [snake=snake,segment amplitude=.2mm,segment length=2mm] (3.7,-1.3)--(6.2,1.3);
\draw [snake=snake,segment amplitude=.2mm,segment length=2mm] (3.7,1.3)--(6.2,-1.3);
\draw [snake=snake,segment amplitude=.2mm,segment length=2mm] (4.3,-0.7)--(4.3,0.7);

\draw (3.4,1.4) node {$\scriptstyle{p'}$}; 
\draw (3.4,-1.4) node {$\scriptstyle{p}$};

\draw (5.3,-1.63) node{$\scriptstyle{a+c}$};
\draw (5.3,1.63) node{$\scriptstyle{a+c}$};

\draw (4.7,0.7) node {$\scriptstyle{r'}$}; 
%\draw (3.4,-1.4) node {$p$};
\draw (4.7,-0.7) node {$\scriptstyle{r'}$};
\draw (4.47,0) node {$\scriptstyle{x}$};
\draw (6.35,0) node {$\scriptstyle{a+c-r'}$};

\draw (6.4,1.4) node {$\scriptstyle{r'}$};
\draw (6.4,-1.4) node {$\scriptstyle{r}$};

\draw (5.3,-0.85) node {$\scriptstyle{\pi_{\gamma/ \varphi}}$};

\draw (3.8,-0.9) node {$\scriptstyle{\pi_{\bar{\varphi}}}$};

\draw (8.5,0.3) node {{\tiny{Thick R2}}};
\draw (8.5,0) node {$=$};
\end{tikzpicture}
\end{center}

\begin{center}
\begin{tikzpicture}

\draw (8.5,0) node {$=\displaystyle{\sum_{\varphi\subset\gamma}\sum_{w\in P(r',x)}}$};

\draw (12,-1.5)--(12,-1);
\draw (12,1)--(12,1.5);
\draw (12,-1)..controls (13,0).. (12,1);
\draw (12,-1)..controls (11,0).. (12,1);
\draw [snake=snake,segment amplitude=.2mm,segment length=2mm] (10,-1.3)--(14,1.3);
\draw [snake=snake,segment amplitude=.2mm,segment length=2mm] (10,1.3)--(14,-1.3);
\draw [snake=snake,segment amplitude=.2mm,segment length=2mm] (10.5,-1)--(10.5,1);

\draw (9.8,1.4) node {$\scriptstyle{p'}$}; 
\draw (9.8,-1.4) node {$\scriptstyle{p}$};

\draw (14.2,1.4) node {$\scriptstyle{r}$}; 
%\draw (3.4,-1.4) node {$p$};
\draw (14.2,-1.4) node {$\scriptstyle{r'}$};
\draw (10.3,-0.6) node {$\scriptstyle{x}$};
\draw (13.3,0) node {$\scriptstyle{a+c-r'}$};
\draw (10.15,0) node {$\scriptstyle{\pi_{\hat{w}}}$};
\draw (10.2,-1.5) node {$\scriptstyle{\pi_{\bar{\varphi}}}$};
\draw (11.7,0.9) node {$\scriptstyle{r'}$};
\draw (12.4,-1.5) node {$\scriptstyle{a+c}$};
\draw (11.6,-0.9) node {$\scriptstyle{\pi_w}$};
\draw (12.6,-1) node {$\scriptstyle{\pi_{\gamma/\varphi}}$};
%\draw (6.4,1.4) node {$r'$};
%\draw (6.4,-1.4) node {$r$};

%\draw (5.3,-1) node {$\pi_{\gamma\over p}$};

%\draw (3.8,-0.8) node {$\pi_{\bar{p}}$};

\filldraw[black]  (11.7,-0.7) circle (1.7pt)
                  (10.5,0) circle (1.7pt)

(10.2,-1.17) circle (1.7pt)

(12.33,-0.7) circle (1.7pt);

\draw (16,0.4) node {{\tiny{Thick R3 + Pitchfork}}};
\draw (16,0) node {$=$};

\end{tikzpicture}
\end{center}

\begin{center}
\begin{tikzpicture}[scale=1.1]
\draw (0,0) node {$=\displaystyle{\sum_{\varphi\subset\gamma}\sum_{w\in P(r',x)}\sum_{i=0}^{r}\sum_{f_1,f_2,f_3\in P(i,r'-i)}} c^{K_i}_{f_1 f_2 f_3}$};
\draw [snake=snake,segment amplitude=.2mm,segment length=2mm] (3,-1.3)--(3.5,-1);
\draw [snake=snake,segment amplitude=.2mm,segment length=2mm] (3.5,-1)--(6.5,1.3);
\draw [snake=snake,segment amplitude=.2mm,segment length=2mm] (3,1.3)--(3.5,1);
\draw [snake=snake,segment amplitude=.2mm,segment length=2mm] (3.5,1)--(6.5,-1.3);
\draw [snake=snake,segment amplitude=.2mm,segment length=2mm] (3.5,-1)--(3.5,1);
\draw [snake=snake,segment amplitude=.2mm,segment length=2mm] (4,-0.6)--(4,0.6);
\draw [snake=snake,segment amplitude=.2mm,segment length=2mm] (6.2,-1.06)--(6.2,1.06);

\draw (5.2,-1.5)--(5.2,-1.25);
\draw (5.2,1.5)--(5.2,1.25);
\draw (5.2,-1.25).. controls (5.42,0)..(5.2,1.25);
\draw (5.2,-1.25).. controls (7,0)..(5.2,1.25);

\filldraw[black]  (3.2,-1.2) circle (1.7pt)
                  (3.5,0) circle (1.7pt)
                  (4,0) circle (1.7pt)
                   
(6.2,0.85) circle (1.7pt)
                   (5.3,-0.7) circle (1.7pt)
                   (5.27,-1) circle (1.7pt)
                   (5.58,-0.97) circle (1.7pt);

\draw (2.8,-1.5) node {$\scriptstyle{p}$};
\draw (2.8,1.5) node {$\scriptstyle{p'}$};

\draw (3.35,0.6) node {$\scriptstyle{x}$};
\draw (3.25,0) node {$\scriptstyle{\pi_{\hat{w}}}$};
\draw (3.5,-1.33) node {$\scriptstyle{\pi_{\bar{\varphi}}}$};

\draw (3.9,0.4) node {$\scriptstyle{i}$};
\draw (4.3,0) node {$\scriptstyle{\pi_{f_1}}$};

\draw (4.5,0.52) node {$\scriptstyle{r-i}$};
\draw (4.5,-0.52) node {$\scriptstyle{r'-i}$};

\draw (5.2,-1.7) node {$\scriptstyle{a+c}$};

\draw (5,-0.7) node {$\scriptstyle{\pi_{\bar{f_3}}}$};
\draw (4.95,-1) node {$\scriptstyle{\pi_w}$};

\draw (5,1) node {$\scriptstyle{r'}$};

\draw (5.81,1.2) node {$\scriptscriptstyle{a+c-r'}$};

%\draw (5.4,-0.2) node {$\scriptstyle{i}$};

\draw (6.5,0.85) node {$\scriptscriptstyle{\pi_{\scriptscriptstyle{f_2}}}$};
\draw (6.05,0) node {$\scriptstyle{i}$};

\draw (5.7,-1.2) node {$\scriptstyle{\pi_{\gamma/\varphi}}$};

\draw (6.7,-1.5) node {$\scriptstyle{r}$};
\draw (6.7,1.5) node {$\scriptstyle{r'}$};

\draw (8,0.4) node {{\tiny{Thick R2}}};
\draw (8,0) node {$=$};

\end{tikzpicture}
\end{center}

\begin{center}
\begin{equation}\label{velsl}
\begin{tikzpicture}[scale=1.1]
\draw (-2,0) node {$=\displaystyle{\sum_{\varphi\subset\gamma}\sum_{w\in P(r',x)}\sum_{i=0}^{r}\sum_{f_1,f_2,f_3\in P(i,r'-i)}\sum_{y\in P(a+c-r',i)}} c^{K_i}_{f_1 f_2 f_3}$};
\draw [snake=snake,segment amplitude=.2mm,segment length=2mm] (3,-1.3)--(3.5,-1);
\draw [snake=snake,segment amplitude=.2mm,segment length=2mm] (3.5,-1)--(6.5,1.3);
\draw [snake=snake,segment amplitude=.2mm,segment length=2mm] (3,1.3)--(3.5,1);
\draw [snake=snake,segment amplitude=.2mm,segment length=2mm] (3.5,1)--(6.5,-1.3);

\draw [snake=snake,segment amplitude=.2mm,segment length=2mm] (3.5,-1)--(3.5,1);

\draw [snake=snake,segment amplitude=.2mm,segment length=2mm] (3.5,-0.7)--(4,-0.4);

\draw [snake=snake,segment amplitude=.2mm,segment length=2mm] (4,-0.4)--(4,0.4);

\draw [snake=snake,segment amplitude=.2mm,segment length=2mm] (3.5,0.7)--(4,0.4);

\draw [snake=snake,segment amplitude=.2mm,segment length=2mm] (6.2,-1.06)--(6.2,1.06);

\draw (5.2,-2)--(5.2,-1.75);
\draw (5.2,2)--(5.2,1.75);

\draw (5.2,-1.75).. controls (5.15,-1.2)..(5.7,-1);
\draw (5.2,-1.75).. controls (5.7,-1.55)..(5.7,-1);

\draw (5.7,-1).. controls (5.8,0.37)..(5.2,1.75);

\filldraw[black]  (3.2,-1.2) circle (1.7pt)
                  (3.5,0) circle (1.7pt)
                  (4,0) circle (1.7pt)
                   
(6.2,0.35) circle (1.7pt)
(6.2,-0.35) circle (1.7pt)

                   (5.3,-1.18) circle (1.7pt)
                   (5.2,-1.55) circle (1.7pt)
                   (5.67,-1.27) circle (1.7pt)

                   (5.5,-1.61) circle (1.7pt);

\draw (5,-1.18) node {$\scriptscriptstyle{\pi_{\scriptscriptstyle{\bar{f_3}}}}$};
\draw (4.9,-1.55) node {$\scriptscriptstyle{\pi_{\scriptscriptstyle{w}}}$};
\draw (5.9,-1.27) node {$\scriptscriptstyle{\pi_{\scriptscriptstyle{y}}}$};
\draw (5.9,-1.7) node {$\scriptscriptstyle{\pi_{\scriptscriptstyle{\gamma/\varphi}}}$};

\draw (5.5,-0.95) node {$\scriptstyle{r'}$};
\draw (6.2,-1) node {$\scriptscriptstyle{a+c-r'}$};

\draw (6.5,0.35) node {$\scriptscriptstyle{\pi_{\scriptscriptstyle{f_2}}}$};
\draw (6.5,-0.35) node {$\scriptscriptstyle{\pi_{\scriptscriptstyle{\hat{y}}}}$};

\draw (2.8,-1.5) node {$\scriptstyle{p}$};
\draw (2.8,1.5) node {$\scriptstyle{p'}$};

\draw (3.35,0.4) node {$\scriptstyle{x}$};
\draw (3.25,0) node {$\scriptstyle{\pi_{\hat{w}}}$};
\draw (3.5,-1.33) node {$\scriptstyle{\pi_{\bar{\varphi}}}$};

\draw (4.1,0.3) node {$\scriptstyle{i}$};
\draw (4.3,0) node {$\scriptstyle{\pi_{f_1}}$};

\draw (4.5,0.52) node {$\scriptstyle{r-i}$};
\draw (4.5,-0.52) node {$\scriptstyle{r'-i}$};

\draw (5.2,-2.25) node {$\scriptstyle{a+c}$};

%\draw (5,-0.7) node {$\scriptstyle{\pi_{\bar{f_3}}}$};
%\draw (4.95,-1) node {$\scriptstyle{\pi_w}$};

%\draw (5,1) node {$\scriptstyle{r'}$};

%\draw (5.81,1.2) node {$\scriptscriptstyle{a+c-r'}$};

%\draw (5.4,-0.2) node {$\scriptstyle{i}$};

\draw (6.35,0.8) node {$\scriptstyle{i}$};

%\draw (5.7,-1.2) node {$\scriptstyle{\pi_{\gamma/\varphi}}$};

\draw (6.7,-1.5) node {$\scriptstyle{r}$};
\draw (6.7,1.5) node {$\scriptstyle{r'}$};
\end{tikzpicture}
\end{equation}
\end{center}
\noindent where $K_i=\underbrace{(r'-i,\ldots,r'-i)}_{i}$, $i>0$, and $K_i=0$ for $i=0$.

Although the last expression has many sums in it, only very few summands can be nonzero. First of all, we have that
\begin{center}
\begin{tikzpicture} [scale=0.55]

\draw (0,-1.5)--(0,1.5);
\draw (-0.5,0.5) node {$\scriptstyle{\pi_w}$};
\draw (-0.5,-0.5) node {$\scriptstyle{\pi_{\bar{f_3}}}$};
\draw (0.2,-1.6) node {$\scriptstyle{r'}$};
\draw (7.2,-1.6) node {$\scriptstyle{r'}$};

\draw (4,0) node {$=\sum\limits_{z\in P(r'\!,x+i)} c_{w,\bar{f_3}}^z$};
\draw (7,-1.5)--(7,1.5);
\draw (7.5,0) node {$\scriptstyle{\pi_z}$};

\filldraw[black]  (0,-0.5) circle (2pt)
 (0,0.5) circle (2pt)
 (7,0) circle (2pt);

 \end{tikzpicture}
\end{center}
 
\noindent and

\begin{center}
\begin{tikzpicture} [scale=0.55]

\draw (0,-1.5)--(0,1.5);
\draw (-0.7,0.5) node {$\scriptstyle{\pi_{\gamma/\varphi}}$};
\draw (-0.5,-0.5) node {$\scriptstyle{\pi_{{y}}}$};
\draw (0.2,-1.6) node {$\scriptstyle{a+c-r'}$};

\draw (4,0) node {$=\sum\limits_{{\scriptsize{\nu\in P(a+c-r',r'-r)}}} c_{\varphi,\nu}^{\gamma}$};
\draw (8,-1.5)--(8,1.5);
\draw (8.5,0.5) node {$\scriptstyle{\pi_\nu}$};
\draw (8.5,-0.5) node {$\scriptstyle{\pi_y}$};
\filldraw[black]  (8,-0.5) circle (2pt)
(0,-0.5) circle (2pt) 
(8,0.5) circle (2pt)
 (0,0.5) circle (2pt);
 
\draw (8.2,-1.6) node {$\scriptstyle{a+c-r'}$};

\draw (16,0) node {$=\sum\limits_{{\scriptstyle{u\in P(a+c-r',r'-r+i)}}}\sum\limits_{{\scriptsize{\nu\in P(a+c-r',r'-r)}}} c_{\varphi,\nu}^{\gamma} c_{y,\nu}^{u}$};
\draw (23,-1.5)--(23,1.5);
\draw (23.5,0) node {$\scriptstyle{\pi_u}$};
%\draw (7.5,-0.5) node {$\scriptstyle{\pi_y}$};
\filldraw[black]  (23,0) circle (2pt);
%(7,-0.5) circle (2pt) 
%(0,0.5) circle (2pt)
% (7,0.5) circle (2pt);
\draw (23.2,-1.6) node {$\scriptstyle{a+c-r'}$};

%\pi_{\gamma/\varphi}\pi_y=\sum\limits_{\nu\in P(a+c-r',r'-r)} c_{\varphi,\nu}^{\gamma} \pi_{\nu}\pi_y = 
%\sum\limits_{u\in P(a+c-r',r'-r+i)}\sum\limits_{\nu\in P(a+c-r',r'-r)} c_{\varphi,\nu}^{\gamma} c_{y,\nu}^{u} \pi_u,$$

 \end{tikzpicture}
\end{center}

\noindent and so the digon on thick line can be written as:
\begin{equation}
%\begin{center}
\begin{tikzpicture} [scale=0.55]
\draw (-1,-1.5)--(-1,-1);
\draw (-1,1.5)--(-1,1);
\draw (-1,-1)..controls (-1.8,0)..(-1,1);
\draw (-1,-1)..controls (-0.2,0)..(-1,1);

\draw (-0.8,-1.6) node {$\scriptstyle{a+c}$};
\draw (-1.2,1.6) node {$\scriptstyle{a+c}$};

\draw (-1.6,0.83) node {$\scriptstyle{r'}$};
\draw (0.1,0.83) node {$\scriptscriptstyle{a+c-r'}$};

\filldraw[black]  (-1.6,0) circle (2pt)
(-1.25,-0.7) circle (2pt)
(-0.4,0) circle (2pt)
(-0.75,-0.7) circle (2pt);

\draw (-2.1,0) node {$\scriptstyle{\pi_{\bar{f_3}}}$};
\draw (-1.7,-0.7) node {$\scriptstyle{\pi_{w}}$};
\draw (0,0) node {$\scriptstyle{\pi_{y}}$};
\draw (0,-0.7) node {$\scriptstyle{\pi_{\gamma/\varphi}}$};

\draw (9,0) node {$\scriptstyle{= \sum\limits_{z\in P(r',x+i)} \sum\limits_{{\scriptsize{u\in P(a+c-r',r'-r+i)}}} \sum\limits_{\nu\in P(a+c-r',r'-r)} c_{\varphi,\nu}^{\gamma} c_{y,\nu}^{u} c_{w,\bar{f_3}}^z}$};
\draw (19,-1.5)--(19,-1);
\draw (19,1.5)--(19,1);
\draw (19,-1)..controls (18.2,0)..(19,1);
\draw (19,-1)..controls (19.8,0)..(19,1);

\draw (19.2,-1.6) node {$\scriptstyle{a+c}$};
\draw (18.8,1.6) node {$\scriptstyle{a+c}$};

\draw (18.4,0.83) node {$\scriptstyle{r'}$};
\draw (20.1,0.83) node {$\scriptscriptstyle{a+c-r'}$};

%\draw (-0.7,-0.7) node {$\scriptstyle{\pi_{w}}$};

\draw (20,0) node {$\scriptstyle{\pi_{u}}$};
\draw (18,0) node {$\scriptstyle{\pi_{z}}$};

\filldraw[black]  (18.4,0) circle (2pt)
%(-0.25,-0.7) circle (2pt)
(19.6,0) circle (2pt);
%(0.25,-0.7) circle (2pt);

 \end{tikzpicture}
%\end{center}
\end{equation}
 %(\pi_{z} \quad,\quad \pi_u) 

%*****************************************************

Now, since $r'\ge r \ge i$ and by assumption  $b\le a+c$, we have that $x+i=p-r'+i\le p-r'+r =b-r'\le a+c-r'$ and $r'-r+i\le r'$, and so $z\in P(r',a+c-r')$ and $u\in P(a+c-r',r')$. Thus, by Lemma \ref{vazlem}, we have that the last diagram can be nonzero only when $|z|+|u|=r'(a+c-r')$, i.e.  we must have 
\begin{equation}\label{vazna1}
|w|+|f_3|+|y|+|\gamma |-| \varphi |=r'(a+c-r').
\end{equation}

As for the digon on curly lines:

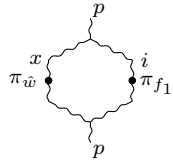
\begin{center}
\begin{tikzpicture} [scale=0.55]
%\draw (3.5,-0.3) node {$\displaystyle \delta_{r,r'}\delta_{\alpha,\alpha'}(-1)^{|\alpha|} \sum\limits_{w\in P(r,x)} \sum\limits_{y\in P(a+c-r,r)}$};
\draw [snake=snake,segment amplitude=.2mm,segment length=2mm](11,-1.5)--(11,-1); 
\draw [snake=snake,segment amplitude=.2mm,segment length=2mm](11,-1)--(10,-0.5);
\draw [snake=snake,segment amplitude=.2mm,segment length=2mm](11,-1)--(12,-0.5);
\draw [snake=snake,segment amplitude=.2mm,segment length=2mm](10,-0.5)--(10,0.5);
\draw [snake=snake,segment amplitude=.2mm,segment length=2mm](12,-0.5)--(12,0.5);
\draw [snake=snake,segment amplitude=.2mm,segment length=2mm](12,0.5)--(11,1);
\draw [snake=snake,segment amplitude=.2mm,segment length=2mm](10,0.5)--(11,1);
\draw [snake=snake,segment amplitude=.2mm,segment length=2mm](11,1.5)--(11,1);

\draw (11.2,-1.7) node {$\scriptstyle{p}$};
\draw (11.2,1.7) node {$\scriptstyle{p}$};
\draw (9.7,0.5) node {$\scriptstyle{x}$};
\draw (12.3,0.5) node {$\scriptstyle{i}$};
\draw (9.4,0) node {$\scriptstyle{\pi_{\hat{w}}}$};
\draw (12.6,-0.1) node {$\scriptstyle{\pi_{f_1}}$};
\filldraw[black]  (10,0) circle (2.2pt)
(12,0) circle (2.2pt);
\end{tikzpicture}
\end{center}

\noindent again by Lemma \ref{vazlem} it can be nonzero only when $|\hat{w}|+|f_1|\ge xi$, i.e.: 
\begin{equation}\label{vazna2}
r'x-|w|+|f_1|\ge xi.
\end{equation}
Thus, from (\ref{vazna1}), (\ref{vazna2}), $|f_1|+|f_3|\le |f_1|+|f_2|+|f_3|=i(r'\!\!-\!i)$ and since $y\in P(a\!+\!c\!-\!r',i)$, $\gamma\in P(a\!+\!c\!-\!r'\!-\!p,r'\!\!-\!r)$ and $x=p\!-\!r'$, we obtain:
\[
r'(a+c-r')\le (p-r')(r'-i)+i(r'-i)+(a+c-r')i+(a+c-r'-p)(r'-r)-|\varphi|.
\]
\noindent The last can be rewritten as
\[
|\varphi|+(a+c-i-p)(r-i)+(r'-i)(r'-r)\le 0.
\]
Since $r'\ge r\ge i$ and $a+c\ge b=p+r\ge p+i$, all terms on the left hand side must be equal to zero, i.e. we must have $r'=r=i$ and $\varphi=0$. Moreover, since $r'=r$, we also have $\gamma=0$, and so by Lemma \ref{vazlem} we have that $\alpha'=\alpha$ and $s_{\gamma}=(-1)^{|\alpha|}$. By replacing all this in (\ref{velsl}), it becomes:
\begin{center}
\begin{tikzpicture} [scale=0.55]
\draw (3.5,-0.3) node {$\displaystyle \delta_{r,r'}\delta_{\alpha,\alpha'}(-1)^{|\alpha|} \sum\limits_{w\in P(r,x)} \sum\limits_{y\in P(a+c-r,r)}$};
\draw [snake=snake,segment amplitude=.2mm,segment length=1mm](11,-1.5)--(11,-1); 
\draw [snake=snake,segment amplitude=.2mm,segment length=1mm](11,-1)--(10,-0.5);
\draw [snake=snake,segment amplitude=.2mm,segment length=1mm](11,-1)--(12,-0.5);
\draw [snake=snake,segment amplitude=.2mm,segment length=1mm](10,-0.5)--(10,0.5);
\draw [snake=snake,segment amplitude=.2mm,segment length=1mm](12,-0.5)--(12,0.5);
\draw [snake=snake,segment amplitude=.2mm,segment length=1mm](12,0.5)--(11,1);
\draw [snake=snake,segment amplitude=.2mm,segment length=1mm](10,0.5)--(11,1);
\draw [snake=snake,segment amplitude=.2mm,segment length=1mm](11,1.5)--(11,1);

\draw (11.2,-1.6) node {$\scriptstyle{p}$};
\draw (11.2,1.6) node {$\scriptstyle{p}$};
\draw (9.8,0.5) node {$\scriptstyle{x}$};
\draw (12.2,0.5) node {$\scriptstyle{r}$};
\draw (9.4,0) node {$\scriptstyle{\pi_{\hat{w}}}$};
\filldraw[black]  (10,0) circle (2.2pt);

\draw (14,-1.5)--(14,-1);
\draw (14,1.5)--(14,1);
\draw (14,-1)..controls (13.2,0)..(14,1);
\draw (14,-1)..controls (14.8,0)..(14,1);

\draw (14.2,-1.68) node {$\scriptstyle{a+c}$};
\draw (14.2,1.68) node {$\scriptstyle{a+c}$};
\draw (13.3,-0.5) node {$\scriptstyle{r}$};
\draw (15.1,-0.6) node {$\scriptscriptstyle{a+c-r}$};
\draw (12.9,0) node {$\scriptstyle{\pi_{{w}}}$};
\draw (15,0) node {$\scriptstyle{\pi_{{y}}}$};
\filldraw[black]  (13.4,0) circle (2pt)
(14.58,0) circle (2pt);

\draw [snake=snake,segment amplitude=.2mm,segment length=1mm](16.5,-1.5)--(16.5,1.5);
\draw (16.7,-1.6) node {$\scriptstyle{r}$};
 \draw (17,0) node {$\scriptstyle{\pi_{\hat{y}}}$};
\filldraw[black]  (16.5,0) circle (2pt);

\end{tikzpicture}
\end{center}

Again, by Lemma \ref{vazlem}, the last is nonzero only when $w=0$ and $y=(\underbrace{r,\ldots,r}_{a+c-r})$ and thus (\ref{velsl}) reduces to

\begin{center}
\begin{tikzpicture} [scale=0.55]
\draw (0,0) node {$\displaystyle  \delta_{r,r'}\delta_{\alpha,\alpha'}(-1)^{|\alpha|} (-1)^{r(a+c-r)}$};
\draw [snake=snake,segment amplitude=.2mm,segment length=2mm](5.5,-1.2)--(5.5,1.2);
\draw (6.5,-1.2)--(6.5,1.2); 
\draw [snake=snake,segment amplitude=.2mm,segment length=2mm](7.5,-1.2)--(7.5,1.2);
\draw (5.7,-1.4) node {$\scriptstyle{p}$};
\draw (6.7,-1.4) node {$\scriptstyle{a+c}$};
\draw (7.7,-1.4) node {$\scriptstyle{r}$};

\end{tikzpicture}
\end{center}

\noindent as wanted.   \kraj

\section{Categorification of the higher quantum Serre relations}\label{sec5}

In this section we give a direct categorification of the higher quantum Serre relations for type $A$. This will be done in the homotopy category of $\dot{\U}^+_3$. The higher quantum Serre relations:
\begin{equation}\label{hqs1}
\sum_{i=0}^a {(-1)^i q^{\pm(a-b-1)i}\,\,\, E_1^{(a-i)}E_2^{(b)}E_1^{(i)}}=0, \quad a>b>0,
\end{equation}
 state that certain alternating sum of monomials in $\u3$ is equal to zero. Thus, it is natural to categorify it by building a complex of objects of $\dot{\U}^+_3$ lifting those monomials which is homotopic to zero.

%Thus, let $i$ and $j$ be elements from the index set such that $i\cdot j=-1$. 
For a categorification of the relation (\ref{hqs}) with the plus sign, our goal is to define a complex of the form 
%\begin{equation}\label{cat1}
%{\scriptstyle{\!\!\!\!\!\!\!\!0\rightarrow \E_1^{(a)}\E_2^{(b)} \rightarrow \E_1^{(a-1)}\E_2^{(b)}\E_1 \{a\!-\!b\!-\!1\} 
%\rightarrow \E_1^{(a-2)}\E_2^{(b)}\E_1^{(2)} \{2(a\!-\!b\!-\!1)\} \rightarrow \cdots \rightarrow \E_2^{(b)}\E_1^{(a)} 
%\{a(a-b-1)\}\rightarrow 0,}}
%\end{equation}
\begin{eqnarray}
&&0\rightarrow \E_1^{(a)}\E_2^{(b)} \rightarrow \E_1^{(a-1)}\E_2^{(b)}\E_1 \{a\!-\!b\!-\!1\} 
\rightarrow \E_1^{(a-2)}\E_2^{(b)}\E_1^{(2)} \{2(a\!-\!b\!-\!1)\} \rightarrow \cdots \nonumber\\
&&\quad\quad \cdots \rightarrow \E_1\E_2^{(b)}\E_1^{(a-1)} \{(a-1)(a-b-1)\}\rightarrow \E_2^{(b)}\E_1^{(a)} \{a(a-b-1)\}\rightarrow 0, \label{cat1}
\end{eqnarray}
that is homotopic to zero.

\begin{theorem}\label{te5}
The following complex 
\begin{center}
\begin{tikzpicture}[scale=0.55]
% \filldraw[black] (-0.67,0.5) circle (2pt)
%				 (0.67,-0.5) circle (2pt);
				 %(1,-3) circle (2pt)
				 %(1,3) circle (2pt)
				 %(-1,0) circle (2pt)
				 %(1,0) circle (2pt);
 
\draw (-19,0) node {$0$};
\draw [->] (-18.5,0)--(-17.5,0);

\draw (-15.8,0) node {$\E_1^{(a)}\E_2^{(b)}$};
 
 \draw [->] (-14,0)--(-12,0);
 \draw (-13.5,0.5)--(-13.5,2);
 \draw (-12.5,1.5)--(-12.5,2);
 \draw (-13.5,0.7)--(-12.5,1.5);
 \draw [snake=snake,segment amplitude=.2mm,segment length=2mm] (-13,0.5)--(-13,2);
 \draw (-13.8,0.4) node {$\scriptscriptstyle a$};
\draw (-13.6,2.2) node {$\scriptscriptstyle a-1$};
\draw (-12.3,2.2) node {$\scriptscriptstyle 1$};
\draw (-12.8,0.4) node {$\scriptscriptstyle b$};

\draw (-8.4,0) node {$\E_1^{(a-1)}\E_2^{(b)}\E_1 \{a\!-\!b\!-\!1\}$};

 \draw [->] (-4.8,0)--(-3,0);
 \draw (-4.5,0.5)--(-4.5,2);
 \draw (-3.5,0.5)--(-3.5,2);
 \draw (-4.5,0.7)--(-3.5,1.5);
 \draw [snake=snake,segment amplitude=.2mm,segment length=2mm] (-4,0.5)--(-4,2);
 \draw (-4.8,0.4) node {$\scriptscriptstyle a-1$};
\draw (-4.6,2.2) node {$\scriptscriptstyle a-2$};
\draw (-3.3,2.2) node {$\scriptscriptstyle 2$};
\draw (-3.3,0.4) node {$\scriptscriptstyle 1$};
\draw (-3.8,2.2) node {$\scriptscriptstyle b$};

\draw (1.2,0) node {$\E_1^{(a-2)}\E_2^{(b)}\E_1^{(2)} \{2(a\!-\!b\!-\!1)\}$};

 \draw [->] (5.5,0)--(6.5,0);
 
 \draw (7.15,0) node {$\cdots$};
 \end{tikzpicture}
\end{center}
\begin{center}
\begin{tikzpicture}[scale=0.55]
\draw (-18.6,0) node {$\cdots$};
\draw [->] (-18,0)--(-17,0);

\draw (-12,0) node {$\E_1\E_2^{(b)}\E_1^{(a-1)} \{(a\!-\!1)(a\!-\!b\!-\!1)\}$};

\draw [->] (-7,0)--(-5,0);
 \draw (-6.5,0.5)--(-6.5,1);
 \draw (-5.5,0.5)--(-5.5,2);
 \draw (-6.5,1)--(-5.5,1.5);
 \draw [snake=snake,segment amplitude=.2mm,segment length=2mm] (-6,0.5)--(-6,2);
 \draw (-6.8,0.4) node {$\scriptscriptstyle 1$};
\draw (-5.3,2.2) node {$\scriptscriptstyle a$};
\draw (-6,2.2) node {$\scriptscriptstyle b$};
\draw (-5.3,0.4) node {$\scriptscriptstyle a-1$};

\draw (-1,0) node {$\E_2^{(b)}\E_1^{(a)} \{a(a\!-\!b\!-\!1)\}$};

\draw [->] (2.5,0)--(3.5,0);
\draw (4.1,0) node {$0$};
\end{tikzpicture}
\end{center}

 is homotopic to zero.
\end{theorem}

\textbf{Proof:}
Denote
%
%
%\begin{center}
%\begin{tikzpicture}[scale=0.5]
%\draw (-18.5,0) node {$C_i:=$};
%\draw (-17.5,0)--(-15.5,0);
%\draw (-17.5,0) node {\it{'small{Ci:=}
%\draw (-16.5,-1.5)--(-16.5,1.5);
%\draw [snake=snake,segment amplitude=.2mm,segment length=2mm] (-15.5,-1.5)--(-15.5,1.5);
%\draw (-14.5,-1.5)--(-14.5,1.5);
%\draw (-16.9,-1.7) node {\it{\small{$\textstyle a\!-\!i$}}};
%\draw (-15.3,-1.7) node {\it{\small{$b$}}};
%\draw (-14.3,-1.7) node {\it{\small{$i$}}};
%\draw (-8.5,0) node {{$\displaystyle{\{i(a\!-\!b\!-\!1)\},\quad\quad i=0,\ldots,a.}$}};
%\end{tikzpicture}
%\end{center}
%i.e. 
$C_i:=\E_1^{(a-i)}\E_2^{(b)}\E_1^{(i)} \{i(a-b-1)\}$,
and 

\begin{center}
\begin{tikzpicture}[scale=0.5]
\draw (-18.5,0) node {$d_i:=$};
%\draw (-17.5,0)--(-15.5,0);
%\draw (-17.5,0) node {\it{'small{Ci:=}
\draw (-16.5,-1.5)--(-16.5,1.5);
\draw [snake=snake,segment amplitude=.2mm,segment length=2mm] (-15,-1.5)--(-15,1.5);
\draw (-13.5,-1.5)--(-13.5,1.5);
\draw (-16.5,-1)--(-13.5,1);
\draw (-16.9,-1.7) node {\it{\small{$\scriptstyle a\!-\!i$}}};
\draw (-16.9,1.7) node {\it{\small{$\scriptstyle a\!-\!(i\!+\!1)$}}};
\draw (-14.8,-1.7) node {\it{\small{$\scriptstyle b$}}};
\draw (-13.8,-1.7) node {\it{\small{$\scriptstyle i$}}};
\draw (-13.8,1.7) node {\it{\small{$\scriptstyle{i+1}$}}};
\draw (-15.8,0) node {\it{\small{$\scriptstyle{1}$}}};
\draw (-7,0) node {$:\quad C_i\longrightarrow C_{i+1},\quad i=0,\ldots,a-1.$};
\end{tikzpicture}
\end{center}
Thus we have to show that the complex
\[
\mathcal{C}:\quad 0 \rightarrow C_0 \xrightarrow{d_0} C_1 \xrightarrow{d_1}  \cdots \xrightarrow{d_{a-1}} C_a \rightarrow 0,  
\]
is homotopic to zero.

First of all, $\mathcal{C}$ is indeed a complex, since

\begin{center}
\begin{tikzpicture}[scale=0.5]
%\draw (-18.5,0) node {$d_i:=$};
%\draw (-17.5,0)--(-15.5,0);
%\draw (-17.5,0) node {\it{'small{Ci:=}
\draw (-21.5,0) node {$d_{i+1}d_i=$};
\draw (-16.5,-1.5)--(-16.5,1.5);
\draw [snake=snake,segment amplitude=.2mm,segment length=2mm] (-15,-1.5)--(-15,1.5);
\draw (-13.5,-1.5)--(-13.5,1.5);

\draw (-16.5,-0.7)--(-13.5,-0.2);
\draw (-16.5,0.7)--(-13.5,1.2);
\draw (-17.3,-1.7) node {\it{\small{$\scriptstyle{a-i}$}}};
\draw (-17.7,0) node {\it{\small{$\scriptstyle{a-(i+1)}$}}};
\draw (-17.7,1.7) node {\it{\small{$\scriptstyle{a-(i+2)}$}}};
\draw (-14.8,-1.7) node {\it{\small{$\scriptstyle{b}$}}};
\draw (-12.8,-1.7) node {\it{\small{$\scriptstyle{i}$}}};
\draw (-12.8,1.7) node {\it{\small{$\scriptstyle{i+2}$}}};
\draw (-12.8,0.1) node {\it{\small{$\scriptstyle{i+1}$}}};
\draw (-15.8,-0.2) node {\it{\small{$\scriptstyle{1}$}}};
\draw (-15.8,1.1) node {\it{\small{$\scriptstyle{1}$}}};

\draw (-11,0) node {$=$};

\draw (-9,-1.5)--(-9,1.5);
\draw [snake=snake,segment amplitude=.2mm,segment length=2mm] (-7.5,-1.5)--(-7.5,1.5);
\draw (-6,-1.5)--(-6,1.5);

\draw (-9,-0.4)--(-8.75,-0.2);
\draw (-8.75,-0.2)..controls (-8, -0.3 ).. (-7.75,0.3);
\draw (-8.75,-0.2)..controls (-8.5,0.3).. (-7.75,0.3);
\draw (-7.75,0.3)--(-6,1);

\draw (-9.6,-1.7) node {\it{\small{$\scriptstyle{a-i}$}}};
%\draw (-9.7,0) node {\it{\small{$\scriptstyle{a-(i+1)}$}}};
\draw (-10,1.7) node {\it{\small{$\scriptstyle{a-(i+2)}$}}};
\draw (-7.15,-1.7) node {\it{\small{$\scriptstyle{b}$}}};
\draw (-5.65,-1.7) node {\it{\small{$\scriptstyle{i}$}}};
\draw (-5.3,1.7) node {\it{\small{$\scriptstyle{i+2}$}}};
%\draw (-5.8,0.1) node {\it{\small{$\scriptstyle{i+1}$}}};
\draw (-8.82,-0.55) node {\it{\small{$\scriptstyle{2}$}}};
\draw (-8.2,-0.5) node {\it{\small{$\scriptstyle{1}$}}};
\draw (-8.5,0.5) node {\it{\small{$\scriptstyle{1}$}}};
\draw (-6.8,1) node {\it{\small{$\scriptstyle{2}$}}};
\draw (-4,0) node {$= \,\,\, 0.$};

\end{tikzpicture}
\end{center}

In order to prove that $\mathcal{C}$ is homotopic to zero, we are left with defining morphisms 
$h_i: C_{i+1} \rightarrow C_i$, for $i=0,\dots,a-1$, such that 
\begin{eqnarray}
h_id_i+d_{i-1}h_{i-1}=\Id\nolimits_{C_i},&\quad i=1,\ldots ,a-1, \label{hom1}\\ 
h_0d_0=\Id\nolimits_{C_0},& \label{hom2}\\
d_{a-1}h_{a-1}=\Id\nolimits_{C_a}.& \label{hom3}
\end{eqnarray}

To that end, for every $i=0,\dots,a-1$, we define $h_i: C_{i+1} \rightarrow C_i$, as follows:

\begin{center}
\begin{tikzpicture}[scale=0.5]
%\draw (-18.5,0) node {$d_i:=$};
%\draw (-17.5,0)--(-15.5,0);
\draw (-20.2,0) node {$h_i:=(-1)^{a-1-i}$};
\draw (-17,-1.5)--(-17,1.5);
\draw [snake=snake,segment amplitude=.2mm,segment length=2mm] (-15.5,-1.5)--(-15.5,1.5);
\draw (-14,-1.5)--(-14,1.5);

\draw (-17,1)--(-14,-0.5);

%\draw (-16.5,0.7)--(-13.5,1.2);
\draw (-17.7,1.7) node {\it{\small{$\scriptstyle{a-i}$}}};
%\draw (-17.7,0) node {\it{\small{$\scriptstyle{a-(i+1)}$}}};
\draw (-18,-1.7) node {\it{\small{$\scriptstyle{a-(i+1)}$}}};

\draw (-15.3,-1.7) node {\it{\small{$\scriptstyle{b}$}}};

\draw (-13.8,1.7) node {\it{\small{$\scriptstyle{i}$}}};
\draw (-13.8,-1.7) node {\it{\small{$\scriptstyle{i+1}$}}};
%\draw (-12.8,0.1) node {\it{\small{$\scriptstyle{i+1}$}}};

\draw (-16.3,0) node {\it{\small{$\scriptstyle{a\!-\!b\!-\!1}$}}};

\draw (-14.5,-0.6) node {\it{\small{$\scriptstyle{1}$}}};

\filldraw[black] (-16,0.5) circle (2pt);

\end{tikzpicture}
\end{center}

Then, for every $i=1,\ldots,a-1$, we have

\begin{center}
\begin{tikzpicture}[scale=0.5]
%\draw (-18.5,0) node {$d_i:=$};
%\draw (-17.5,0)--(-15.5,0);
\draw (-22,0) node {$h_id_i+d_{i-1}h_{i-1}=(-1)^{a-1-i}$};
\draw (-17,-1.5)--(-17,1.5);
\draw [snake=snake,segment amplitude=.2mm,segment length=2mm] (-15.5,-1.5)--(-15.5,1.5);
\draw (-14,-1.5)--(-14,1.5);

\draw (-17,1)--(-14,0.2);
\draw (-17,-1)--(-14,-0.2);

%\draw (-16.5,0.7)--(-13.5,1.2);
\draw (-17.7,1.7) node {\it{\small{$\scriptstyle{a-i}$}}};
%\draw (-17.7,0) node {\it{\small{$\scriptstyle{a-(i+1)}$}}};
\draw (-18,-1.7) node {\it{\small{$\scriptstyle{a-i}$}}};

\draw (-15.3,-1.7) node {\it{\small{$\scriptstyle{b}$}}};

\draw (-13.8,1.7) node {\it{\small{$\scriptstyle{i}$}}};
\draw (-13.8,-1.7) node {\it{\small{$\scriptstyle{i}$}}};
\draw (-13.5,0) node {\it{\small{$\scriptstyle{i\!+\!1}$}}};

\filldraw[black] (-16.2,0.77) circle (2pt);

\draw (-16.16,0.4) node {\it{\small{$\scriptscriptstyle{a\!-\!b\!-\!1}$}}};

\draw (-14.5,-0.6) node {\it{\small{$\scriptstyle{1}$}}};
\draw (-14.5,0.6) node {\it{\small{$\scriptstyle{1}$}}};

\draw (-11,0) node {$+ (-1)^{a-i}$};

\draw (-9,-1.5)--(-9,1.5);
\draw [snake=snake,segment amplitude=.2mm,segment length=2mm] (-7.5,-1.5)--(-7.5,1.5);
\draw (-6,-1.5)--(-6,1.5);

\draw (-9,-0.2)--(-6,-1);
\draw (-9,0.2)--(-6,1);

%\draw (-16.5,0.7)--(-13.5,1.2);

\draw (-9.6,1.7) node {\it{\small{$\scriptstyle{a-i}$}}};
%\draw (-17.7,0) node {\it{\small{$\scriptstyle{a-(i+1)}$}}};
\draw (-9.6,-1.7) node {\it{\small{$\scriptstyle{a-i}$}}};

\draw (-7.3,-1.7) node {\it{\small{$\scriptstyle{b}$}}};

\draw (-6.7,1.2) node {\small $\scriptstyle{1}$};
\draw (-6.7,-1.2) node {\small $\scriptstyle{1}$};

\draw (-5.8,1.7) node {\it{\small{$\scriptstyle{i}$}}};
\draw (-5.8,-1.7) node {\it{\small{$\scriptstyle{i}$}}};
\draw (-5.5,0.1) node {\it{\small{$\scriptstyle{i\!-\!1}$}}};

\draw (-8.2,-1) node {\it{\small{$\scriptscriptstyle{a\!-\!b\!-\!1}$}}};

%\draw (-14.5,-0.6) node {\it{\small{$\scriptstyle{1}$}}};
%\draw (-14.5,0.6) node {\it{\small{$\scriptstyle{1}$}}};
\filldraw[black] (-8,-0.48) circle (2pt);

\draw (-4,0) node {$=$};

\end{tikzpicture}
\end{center}

\begin{center}
\begin{tikzpicture}[scale=0.5]
%\draw (-18.5,0) node {$d_i:=$};
%\draw (-17.5,0)--(-15.5,0);
\draw (-20,0) node {$=(-1)^{a-1-i}$};
\draw (-17,-1.5)--(-17,1.5);
\draw [snake=snake,segment amplitude=.2mm,segment length=2mm] (-15.8,-1.5)--(-15.8,1.5);
\draw (-14,-1.5)--(-14,1.5);

\draw (-17,0.7)--(-14,-0.7);
\draw (-17,-0.7)--(-14,0.7);

%\draw (-16.5,0.7)--(-13.5,1.2);
\draw (-17.7,1.7) node {\it{\small{$\scriptstyle{a-i}$}}};
%\draw (-17.7,0) node {\it{\small{$\scriptstyle{a-(i+1)}$}}};
\draw (-18,-1.7) node {\it{\small{$\scriptstyle{a-i}$}}};

\draw (-15.6,-1.7) node {\it{\small{$\scriptstyle{b}$}}};

\draw (-13.8,1.7) node {\it{\small{$\scriptstyle{i}$}}};
\draw (-13.8,-1.7) node {\it{\small{$\scriptstyle{i}$}}};
\draw (-13.5,0) node {\it{\small{$\scriptstyle{i\!-\!1}$}}};

\draw (-16.4,0.9) node {\it{\small{$\scriptscriptstyle{a\!-\!b\!-\!1}$}}};

\draw (-14.5,-0.8) node {\it{\small{$\scriptstyle{1}$}}};
\draw (-14.5,0.8) node {\it{\small{$\scriptstyle{1}$}}};

\filldraw[black] (-16.5,0.48) circle (2pt);

\draw (-11,0) node {$+ (-1)^{a-i}$};

\draw (-9,-1.5)--(-9,1.5);
\draw [snake=snake,segment amplitude=.2mm,segment length=2mm] (-7.2,-1.5)--(-7.2,1.5);
\draw (-6,-1.5)--(-6,1.5);

\draw (-9,-0.7)--(-6,0.7);
\draw (-9,0.7)--(-6,-0.7);

%\draw (-16.5,0.7)--(-13.5,1.2);
\draw (-9.6,1.7) node {\it{\small{$\scriptstyle{a-i}$}}};
%\draw (-17.7,0) node {\it{\small{$\scriptstyle{a-(i+1)}$}}};
\draw (-9.6,-1.7) node {\it{\small{$\scriptstyle{a-i}$}}};

\draw (-7,-1.7) node {\it{\small{$\scriptstyle{b}$}}};

\draw (-5.8,1.7) node {\it{\small{$\scriptstyle{i}$}}};
\draw (-5.8,-1.7) node {\it{\small{$\scriptstyle{i}$}}};
\draw (-5.5,0.1) node {\it{\small{$\scriptstyle{i\!-\!1}$}}};

\draw (-6.5,-1) node {\it{\small{$\scriptscriptstyle{a\!-\!b\!-\!1}$}}};

\draw (-8.5,-0.82) node {\it{\small{$\scriptstyle{1}$}}};
\draw (-8.5,0.82) node {\it{\small{$\scriptstyle{1}$}}};

\filldraw[black] (-6.5,-0.45) circle (2pt);

\draw (-4,0) node {$=$};

\end{tikzpicture}
\end{center}
\begin{center}
\begin{tikzpicture}[scale=0.5]

\draw (-12,0) node {$=\quad (-1)^{a-1-i}$};

\draw (-9,-1.5)--(-9,1.5);
\draw [snake=snake,segment amplitude=.2mm,segment length=2mm] (-7.2,-1.5)--(-7.2,1.5);
\draw (-6,-1.5)--(-6,1.5);

\draw (-9,-0.7)--(-6,0.7);
\draw (-9,0.7)--(-6,-0.7);

%\draw (-16.5,0.7)--(-13.5,1.2);
\draw (-9.6,1.7) node {\it{\small{$\scriptstyle{a-i}$}}};
%\draw (-17.7,0) node {\it{\small{$\scriptstyle{a-(i+1)}$}}};
\draw (-9.6,-1.7) node {\it{\small{$\scriptstyle{a-i}$}}};

\draw (-7,-1.7) node {\it{\small{$\scriptstyle{b}$}}};

\draw (-5.8,1.7) node {\it{\small{$\scriptstyle{i}$}}};
\draw (-5.8,-1.7) node {\it{\small{$\scriptstyle{i}$}}};
\draw (-5.5,0.1) node {\it{\small{$\scriptstyle{i\!-\!1}$}}};

\draw (-8.4,1) node {\it{\small{$\scriptscriptstyle{a\!-\!b\!-\!1}$}}};

\draw (-6.5,-1) node {\it{\small{$\scriptstyle{1}$}}};
\draw (-6.5,1) node {\it{\small{$\scriptstyle{1}$}}};

\filldraw[black] (-8.5,0.45) circle (2pt);

\draw (-1,0) node {$\displaystyle{+ (-1)^{a-1-i}\!\!\!\sum_{x+y+z=b-1}}$};

\draw (3,-1.5)--(3,1.5);
\draw [snake=snake,segment amplitude=.2mm,segment length=2mm] (4.5,-1.5)--(4.5,1.5);
\draw (6,-1.5)--(6,1.5);

\draw (3,-0.7)..controls (3.7,-0.7) and (3.7,0.7)..(3,0.7);
\draw (6,-0.7)..controls (5.3,-0.7) and (5.3,0.7)..(6,0.7);

%\draw (-16.5,0.7)--(-13.5,1.2);
\draw (2.4,1.7) node {\it{\small{$\scriptstyle{a-i}$}}};
%\draw (-17.7,0) node {\it{\small{$\scriptstyle{a-(i+1)}$}}};
\draw (2.4,-1.7) node {\it{\small{$\scriptstyle{a-i}$}}};

\draw (4.8,-1.7) node {\it{\small{$\scriptstyle{b}$}}};

\draw (6.2,1.7) node {\it{\small{$\scriptstyle{i}$}}};
\draw (6.2,-1.7) node {\it{\small{$\scriptstyle{i}$}}};
\draw (6.45,0) node {\it{\small{$\scriptstyle{i\!-\!1}$}}};

\draw (3.6,-0.9) node {\it{\small{$\scriptscriptstyle{a\!-\!b\!-\!1\!+\!x}$}}};

\draw (3.5,0.9) node {\it{\small{$\scriptstyle{1}$}}};

\draw (5.5,-0.9) node {\it{\small{$\scriptstyle{1}$}}};

\draw (5.5,0.7) node {\it{\small{$\scriptstyle{y}$}}};

\draw (5,1.3) node {\it{\small{$\scriptstyle{\epsilon_z}$}}};

\filldraw[black] (3.4,-0.48) circle (2pt)
                  (5.57,0.38) circle (2pt)
                  (4.5,1.3) circle (2pt);
 
\draw (8,0) node {$+$};                 
\end{tikzpicture}
\end{center}

\begin{equation}\label{proba1}
\begin{tikzpicture}[scale=0.5]

\draw (-11,0) node {$+ (-1)^{a-i}$};

\draw (-9,-1.5)--(-9,1.5);
\draw [snake=snake,segment amplitude=.2mm,segment length=2mm] (-7.2,-1.5)--(-7.2,1.5);
\draw (-6,-1.5)--(-6,1.5);

\draw (-9,-0.7)--(-6,0.7);
\draw (-9,0.7)--(-6,-0.7);

%\draw (-16.5,0.7)--(-13.5,1.2);
\draw (-9.6,1.7) node {\it{\small{$\scriptstyle{a-i}$}}};
%\draw (-17.7,0) node {\it{\small{$\scriptstyle{a-(i+1)}$}}};
\draw (-9.6,-1.7) node {\it{\small{$\scriptstyle{a-i}$}}};

\draw (-7,-1.7) node {\it{\small{$\scriptstyle{b}$}}};

\draw (-5.8,1.7) node {\it{\small{$\scriptstyle{i}$}}};
\draw (-5.8,-1.7) node {\it{\small{$\scriptstyle{i}$}}};
\draw (-5.5,0.1) node {\it{\small{$\scriptstyle{i\!-\!1}$}}};

\draw (-6.5,-1) node {\it{\small{$\scriptscriptstyle{a\!-\!b\!-\!1}$}}};

\draw (-8.5,-1) node {\it{\small{$\scriptstyle{1}$}}};
\draw (-8.5,1) node {\it{\small{$\scriptstyle{1}$}}};

\filldraw[black] (-6.5,-0.45) circle (2pt);

\end{tikzpicture}
\end{equation}

Here in the second equality we have used Opening of a Thick Edge, while in the third equality we have applied the Thick R3 move on the first term. Now, by Dot Migration, the sum of the first and the third summand from above is equal to:

\begin{equation}\label{proba2}
\begin{tikzpicture}[scale=0.5]

\draw (-1,0) node {$\displaystyle{(-1)^{a-1-i}\!\!\!\sum_{r+s=a-b-2}}$};

\draw (3,-1.5)--(3,1.5);
\draw [snake=snake,segment amplitude=.2mm,segment length=2mm] (5.38,-1.5)--(5.38,1.5);
\draw (6,-1.5)--(6,1.5);

\draw (3,-1)..controls (4.4,-0.5) and (4.4,0.5)..(3,1);
\draw (6,-1)..controls (4.6,-0.5) and (4.6,0.5)..(6,1);

%\draw (-16.5,0.7)--(-13.5,1.2);
\draw (2.4,1.7) node {\it{\small{$\scriptstyle{a-i}$}}};
%\draw (-17.7,0) node {\it{\small{$\scriptstyle{a-(i+1)}$}}};
\draw (2.4,-1.7) node {\it{\small{$\scriptstyle{a-i}$}}};

\draw (5.5,-1.7) node {\it{\small{$\scriptstyle{b}$}}};

\draw (6.2,1.7) node {\it{\small{$\scriptstyle{i}$}}};
\draw (6.2,-1.7) node {\it{\small{$\scriptstyle{i}$}}};
%\draw (6.45,0) node {\it{\small{$\scriptstyle{i\!-\!1}$}}};

%\draw (3.6,-0.9) node {\it{\small{$\scriptscriptstyle{a\!-\!b\!-\!1\!+\!x}$}}};

\draw (4.1,-0.7) node {\it{\small{$\scriptstyle{1}$}}};

\draw (5,-0.7) node {\it{\small{$\scriptstyle{1}$}}};

\draw (4.1,0.7) node {\it{\small{$\scriptstyle{r}$}}};

\draw (5,0.7) node {\it{\small{$\scriptstyle{s}$}}};

\filldraw[black] (3.92,0.4) circle (2pt)
                  (5.12,0.4) circle (2pt);
                  %(4.5,1.3) circle (2pt);
 
\draw (10.7,0) node {$\displaystyle{= (-1)^{a-1-i}\sum_{\scriptscriptstyle{r\!+\!s=a\!-\!b\!-\!2}}\sum_{x=0}^b}$};

\draw (15,-1.5)--(15,1.5);
\draw [snake=snake,segment amplitude=.2mm,segment length=2mm] (16.5,-1.5)--(16.5,1.5);
\draw (18,-1.5)--(18,1.5);

\draw (15,-1)..controls (16.4,-0.5) and (16.4,0.5)..(15,1);
\draw (18,-1)..controls (16.6,-0.5) and (16.6,0.5)..(18,1);

%\draw (-16.5,0.7)--(-13.5,1.2);
\draw (14.4,1.7) node {\it{\small{$\scriptstyle{a-i}$}}};
%\draw (-17.7,0) node {\it{\small{$\scriptstyle{a-(i+1)}$}}};
\draw (14.4,-1.7) node {\it{\small{$\scriptstyle{a-i}$}}};

\draw (16.8,-1.7) node {\it{\small{$\scriptstyle{b}$}}};

\draw (18.2,1.7) node {\it{\small{$\scriptstyle{i}$}}};
\draw (18.2,-1.7) node {\it{\small{$\scriptstyle{i}$}}};
%\draw (6.45,0) node {\it{\small{$\scriptstyle{i\!-\!1}$}}};

%\draw (3.6,-0.9) node {\it{\small{$\scriptscriptstyle{a\!-\!b\!-\!1\!+\!x}$}}};

\draw (16.1,-0.7) node {\it{\small{$\scriptstyle{1}$}}};

\draw (17,-0.7) node {\it{\small{$\scriptstyle{1}$}}};

\draw (16.1,0.7) node {\it{\small{$\scriptstyle{r}$}}};

\draw (17,0.7) node {\it{\small{$\scriptstyle{x\!+\!s}$}}};

\draw (17.15,1.2) node {\it{\small{$\scriptstyle{\epsilon_{b\!-\!x}}$}}};

\filldraw[black] (15.92,0.4) circle (2pt)
                  (17.12,0.4) circle (2pt)
                  (16.5,1.2) circle (2pt);

\end{tikzpicture}
\end{equation}
where we have used the Thick R2 move. 

The last double sum is nonzero only when $r\ge a-i-1$ and $x+s\ge i-1$. Since $x\le b$ and $r+s=a-b-2$, the last implies
that we must have equalities, i.e. that $r=a-i-1$, $x=b$ and $s=i-b-1$. Finally, the last is possible only when $i\ge b+1$, and so (\ref{proba2}) is equal to $\delta_{\{i\ge b+1\}} \Id_{C_i}$.

Analogously, the remaining second term from (\ref{proba1}) is nonzero, only when $a-b-1+x\ge a-i-1$, $y\ge i-1$ and $z\ge 0$, and since $x+y+z=a-b-2$, again we must have all equalities, and in that case the value of the second term is equal to $\Id_{C_i}$. Since $x\ge 0$, we must have $i\le b$, and so altogether we have:
\[
h_id_i+d_{i-1}h_{i-1}=\delta_{\{i\ge b+1\}}\Id\nolimits_{C_i}+\delta_{\{i\le b\}}\Id\nolimits_{C_i}=\Id\nolimits_{C_i},
\]
thus proving (\ref{hom1}). The equalities (\ref{hom2}) and (\ref{hom3}) can be obtained completely analogously. Hence, $\mathcal{C}$ is homotopic to zero, as wanted. \kraj\\

As for the categorification of the higher quantum Serre relations (\ref{hqs1}) with the minus sign, we obtain it by defining a 
complex of the form (\ref{cat1}), but with the arrows pointing in the opposite direction. By exchanging the roles of $d_i$'s and 
$h_i$'s from the Theorem above, we have: 
\begin{theorem}
The following complex
\begin{center}
$
\!\!\!\!\!\!\!\!0\leftarrow \E_1^{(a)}\E_2^{(b)} \xleftarrow{h_1} \E_1^{(a-1)}\E_2^{(b)}\E_1 \{-(a-b-1)\} \xleftarrow{h_2}
\E_1^{(a-2)}\E_2^{(b)}\E_1^{(2)} \{-2(a-b-1)\} \xleftarrow{h_3} \cdots \xleftarrow{h_{a-1}} \E_2^{(b)}\E_1^{(a)} \{-a(a-b-1)\}\leftarrow 0,
$
\end{center}
is homotopic to zero.
\end{theorem}

\vskip 0.3cm
E-mail: {\tt{mstosic@math.ist.utl.pt}} 

Instituto de Sistemas e Rob\'otica and CAMGSD, 
Instituto Superior T\'ecnico, Torre Norte, Piso 7, 
Av. Rovisco Pais, 1049-001 Lisbon, Portugal,

{\it{and}}

Mathematical Institute SANU, Knez Mihailova 36, 11000 Beograd, Serbia.

\end{document}